\documentclass[12pt, draft]{article}
\usepackage{amsmath}
\usepackage{amssymb}
\usepackage{latexsym}
\usepackage{color}
\allowdisplaybreaks[4]
\renewcommand{\theequation}{\arabic{section}.\arabic{equation}}
\newtheorem{Thm}{Theorem}[section]
\newtheorem{Prop}[Thm]{Proposition}
\newtheorem{Lemma}[Thm]{Lemma}

\newtheorem{Def}[Thm]{Definition}
\newtheorem{Remark}[Thm]{Remark}
\newtheorem{rei}[Thm]{Example}

\newcommand{\R}{{\mathbb R}}

\renewcommand{\v}[1]{\vert #1 \vert}
\newcommand{\V}[1]{\Vert #1 \Vert}

\title{Asymptotic behavior of the indicator function in the inverse problem of the wave equation for media with multiple types of cavities}
\author{
Mishio Kawashita\thanks{Partly supported by JSPS KAKENHI Grant Number JP23K03184.}
\thanks{kawasita@hiroshima-u.ac.jp}
\\
\and 
Wakako Kawashita\thanks{Partly supported by JSPS KAKENHI Grant Number JP20K03684.}
\thanks{wakawa@hiroshima-u.ac.jp}
}
\begin{document}
\maketitle
\begin{abstract}
\noindent
{\bf Abstract.} 
In this paper, the inverse problem of the wave equation by the enclosure method for a medium with multiple types of cavities is discussed. In the case considered here, the sign of the indicator function of the enclosure method is not determined and sign cancellation may occur, resulting in loss of information. By examining the top terms of the indicator function in detail, we show that the shortest distance to the cavities can be obtained even in such a case.
\par
\end{abstract}
\vskip 3truepc

\setcounter{equation}{0}
\section{Introduction}
\label{Introduction}

The enclosure method was introduced by M. Ikehata as a nondestructive evaluation method in the elliptic boundary value problems(\cite{Ikehata M 1999, Ikehata M 2000}). 
The study of enclosure method for time dependent problems is started from one-dimensional heat equation \cite{I4}, and extended to the three-dimensional case by \cite{Ikehata-Kawashita1, Ikehata-Kawashita2, Ikehata-Kawashita3}. Ikehata has conducted many other studies on reconstruction for various equations by the enclosure method(see e.g.\cite{Ikehata M rev}). Also for the wave equation, each case of Dirichlet boundary condition, Neumann type boundary condition, dissipative boundary condition and inclusions is investigated in \cite{IE, IEO2, IEO3 bistatic}.
In the enclosure method, it is important to introduce an appropriate indicator function with a parameter for each problem, and the internal information of the object can be extracted from the analysis of the indicator function for large parameter. 
It is known from previous studies that the\lq\lq shortest length'' for each problem can be extracted from the asymptotic behavior of the indicator function
(see.\,e.g. \cite{IE, IEO2, IEO3 bistatic, Ikehata-Kawashita1, Ikehata-Kawashita2, Ikehata-Kawashita3, transmission No1, transmission No2}). 
\par
If we consider the indicator function for the wave equation in a medium with multiple cavities, but all the cavities have only boundaries that satisfy the Dirichlet (resp. Neumann) condition, the sign of the indicator function will be negative (resp. positive) when this parameter tends to $\infty$(cf. \cite{IEO2, IEO3 bistatic}). In these cases, we say that the\lq\lq monotonicity condition'' is said to be satisfied, and information about cavities such as the shortest length can be obtained from the indicator function relatively easily.
On the other hand, if the cavities have multiple types of boundary conditions, for example, if some of the boundaries are Dirichlet conditions and some are Neumann conditions, it is simply not clear whether the indicator function has a definite sign when the parameter is large enough. Thus, if the monotonicity condition is not satisfied, sign cancellation may occur, and internal information may not be extracted from the indicator function. 
For this reason, many studies have been considered under the monotonicity condition. 
However, making such an assumption is nothing but giving a priori information on the subject. \\
\indent
Therefore, the purpose of this paper is to extract information such as the shortest distance from the observation site to the cavities by studying the asymptotic behavior of the indicator function in the case of cavities with multiple types of boundary conditions. This is a typical case where the monotonicity condition is not satisfied.
%
%
%
\par
We consider a medium with cavities consisting of
disjoint union of bounded open sets 
$D_{j}^{\alpha} $ ($j = 1, 2, \ldots, N^\alpha$) 
$\alpha \in \{n_-, n_+, d \}$ with $C^2$ boundaries.  
We put $\Omega = \R^3\setminus\overline{D}$, 
$D = \cup_{\alpha \in \{n_+, n_-, d\}}D^\alpha$ and $D^\alpha = \cup_{j = 1}^{N^\alpha}D_{j}^{\alpha}$ for $\alpha \in \{ n_-, n_+, d\}$, $\partial\Omega = \partial{D}$ 
and $D^n=D^{n_+}\cup D^{n_-}$. 
Since $D_{j}^{\alpha}$ are disjoint each other, the total number of cavities $N^t$ is given by $N^t=N^{n_+}+N^{n_-} +N^d$. 
We also define a boundary operator 
$\partial_{\nu_{x}}u 
= \sum_{k = 1}^3(\nu_{x})_k \partial_{x_k}u$, where
$\nu_{x} = {}^t((\nu_{x})_1, (\nu_{x})_2, (\nu_{x})_3)$ is 
the unit outer normal of $\partial{D}$ from the $D$-side and 
$\partial_{x_k}=\frac{\partial}{\partial x_k}$.  
\par
Take $T > 0$ as a observation time, and consider the following problem: 
\begin{equation}
\left\{\begin{array}{ll}
(\partial_t^2-\gamma_0 \Delta)u(t, x) = 0 &
\qquad\text{ in } (0, T)\times\Omega, \\
(\gamma_0\partial_{\nu_x}-\lambda_1(x)\partial_t-\lambda_0(x))u(t, x) = 0 
& \qquad\text{ on } (0, T)\times\partial{D^{n}}, \\
u(t, x) = 0 
& \qquad\text{ on } (0, T)\times\partial{D^{d}},
\\
u(0, x) = 0, \enskip \partial_tu(0, x) = f(x) & \qquad\text{ on } 
\Omega, 
\end{array}
\right.
\label{time dependent problem with cavities}
\end{equation}
where $\gamma_0 >0$ is a constant and $\lambda_j \in L^\infty(\partial{D^n})$ $(j = 0, 1)$.    
Thus, on $\partial{D^{n}}$  (resp. $\partial{D^{d}}$), 
the Neumann type (resp. Dirichlet) 
boundary conditions are imposed. 
In this article, we assume that the time derivative term 
$\lambda_1(x)\partial_t$ in the boundary condition on $(0,T) \times \partial D^n$
stands for a dissipation effect: 
$$
\lambda_1(x) \geq 0, \enskip x \in \partial{D^n},
$$
which ensure that the total energy does not increase as time goes by.
Assume that $D^{n_+}$ and $D^{n_-}$ have the following properties:
there exists a constant $\mu_1 > 0$ such that
\begin{equation}
\left\{
\begin{array}{l}
\text{$0 \le \lambda_1(x)  < \sqrt{\gamma_0} -\mu_1$ a.e. on $\partial{D}^{n_+}$}, \\
\text{$0 < \sqrt{\gamma_0} + \mu_1 < \lambda_1(x)$ a.e. on $\partial{D}^{n_-}$.}
\end{array}
\right.
\nonumber
\end{equation}
Take an open set $B$ with 
$\overline{B} \subset \Omega$, and put 
$f \in L^2(\R^3)\cap L^{\infty}(\R^3)$ satisfying the emission condition from $B$, which is given by
\begin{equation}
\left\{
\begin{array}{ll}
\text{$f \in C^1(\overline{B})$ with ${\rm supp } f \subset \overline{B} $ and there exists a constant $c > 0$ such that } 
 \\
\text{$f(x) \geq c$ ($x \in \overline{B}$) or $-f(x) \geq c$ ($x \in \overline{B}$).}
\end{array}
\right.
\label{emission condition}
\end{equation}
This condition (\ref{emission condition}) ensures that the waves are emitted from the entire $B$, including up to its boundary. For this initial data $f$, we measure the waves $u(t, x)$ on $B$ from $0$ to $T$. The inverse problem considered here is to obtain information about the cavity $D$ from this measurement. 
\par
We denote $H^1_{0, \partial D^d}(\Omega)$ by the Sobolev space 
consisting of all functions 
$\varphi \in H^1(\Omega)$ with $\varphi = 0$ on 
$\partial{D^{d}}$ in the trace sense. 
It is well known that for any $f \in L^2(\Omega)$, there exists a unique weak solution $u \in L^2(0, T ; H^1_{0, \partial D^d}(\Omega))$ of (\ref{time dependent problem with cavities}) with $\partial_tu \in L^2(0, T ; L^2(\Omega))$ and 
$\partial_t^2u \in L^2(0, T ; (H^1(\Omega))')$ satisfying 
\begin{align*}
\langle &\partial_t^2u(t, \cdot), \varphi \rangle 
+ \int_{\Omega}\gamma_0 \nabla_xu(t, x)\cdot
\nabla_x\varphi(x)dx \\
&+ \frac{d}{dt}\int_{\partial{D}^h}\lambda_1(x)u(t, x)\varphi(x)dS_x + \int_{\partial{D}^h}\lambda_0(x)u(t, x)\varphi(x)dS_x = 0 
\qquad \text{a.e. $t \in (0, T)$}
\end{align*}
for all $\varphi \in H^1_{0, \partial D^d}(\Omega)$, and the initial condition in usual sense (cf. for example \cite{DL}). 
In what follows, we consider these weak solutions as the class of solutions.  
For the background equation, we consider  
\begin{equation}
\begin{array}{ll}
\displaystyle
(\gamma_0\Delta-\tau^2)v(x; \tau) + f(x) = 0 & 
\displaystyle
 \qquad\text{in}\,\,\,\R^3,
\end{array}
\label{reduced background equation}
\end{equation}
and we use the weak solution $v \in H^1(\R^3)$ 
with the kernel representation:
\begin{equation}
v(x; \tau) = \int_{\Omega}\Phi_\tau(x, y)f(y)dy \quad\text{with}\quad
\Phi_\tau(x, y) 
= \frac{1}{4\pi\gamma_0}\frac{e^{-{\tau}\v{x - y}/\sqrt{\gamma_0}}}{\v{x - y}}.
\label{free kernel_rep}
\end{equation}
\indent
As in the usual approach of the enclosure method, 
we introduce the indicator function $I_{\tau}$ defined by
\begin{equation}
I_{\tau} = \int_{\Omega}f(x)({\mathcal L}_Tu(x; \tau) - v(x; \tau))dx, 
\label{Indicator function for Cave&Inclusion}
\end{equation}
where 
\begin{equation}
{\mathcal L}_Tu(x; \tau)= \int_0^Te^{-{\tau}t}u(t, x)dt \qquad(x \in \Omega). 
\label{Image of partial LT for u}
\end{equation}
Note that $I_{\tau}$ is obtained from the measurement $u(t, x)$ for 
$0 \leq t \leq T$ and $x \in B$. 
\par

\begin{Remark} 
$I_{\tau}$ is actually determined only by $(u, \partial_{\nu_x}u)$ on $(0, T)\times\partial{B}$. In fact, it can be seen from the following formula obtained by integration by parts: 
$$
I_{\tau} = \gamma_0\int_{\partial{B}}(v(x; \tau)\partial_{\nu_x}{\mathcal L}_Tu(x; \tau)-\partial_{\nu_x}v(x; \tau){\mathcal L}_Tu(x; \tau) )dS_x.
$$
\end{Remark}
\par
For the problem (\ref{time dependent problem with cavities}), the \lq\lq shortest length'' is defined as the shortest distance from the observation site to the cavities :
\begin{equation}
\left\{
\begin{array}{ll}
&{\displaystyle l_0 ={\rm dist}(D, B) 
= \hskip-10pt \inf_{x \in \overline{D}, y\in \overline{B}}L_0(x,y), } \enskip
\text{where $L_0(x,y)=\v{x - y}$,}
\\
&l_0^{\alpha} ={\rm dist}(D^\alpha, B) 
\text{ for $\alpha \in \{n_+, n_-, d \}$},\quad
l_0^+=l_0^{n_+}, \enskip l_0^-= \min\{ l_0^{n_-}, l_0^d \}.
\end{array}
\right.
\label{shortest lenghs1}
\end{equation}
Hereafter in this article, we call $D^{n_+}$ is positive cavity, $D^{n_-}$ and $D^d$ are negative cavities.     
Using the elliptic estimate method, which has been Ikehata's primary method for analyzing indicator functions, the authors obtained the following results: 
\begin{Thm}\label{mixed but separated case}(\cite{M. and W. Kawashita separated, M. and W. Kawashita combined})
If all of $D^{n_\pm}$, $D^{d}$ and $B$ have $C^1$ boundaries, 
the indicator function (\ref{Indicator function for Cave&Inclusion}) satisfies the following properties:
\par\noindent
(1) For $T < 2l_0/\sqrt{\gamma_0}$, $\lim_{\tau \to \infty}e^{{\tau}T}I_\tau = 0$.
\par\noindent
(2) For $T > 2l_0/\sqrt{\gamma_0}$, 
if (\ref{emission condition}) holds, then we have: 
$$\text{$(ii)_{+}$} \enskip
\lim_{\tau \to \infty}e^{{\tau}T}I_\tau = +\infty \enskip\text{if $l_0^+ =l_0 < l_0^-$,} 
\quad \text{$(ii)_{-}$} \enskip
\lim_{\tau \to \infty}e^{{\tau}T}I_\tau =-\infty \enskip \text{if $l_0^-=l_0 < l_0^+$.}
$$
Further, in each of the cases $(ii)_+$ and $(ii)_{-}$, we also have
\begin{equation*}
\lim_{\tau \to \infty}\frac{\sqrt{\gamma_0}}{2\tau}\log\v{I_\tau} = 
\begin{cases} -l_0^+ & \text{for the case of $(ii)_+$},
\\
-l_0^- & \text{for the case of $(ii)_-$}.
\end{cases}
\end{equation*}
\end{Thm}
\par
In \cite{M. and W. Kawashita separated}, the authors deal with the case where not only cavities but also heterogeneous parts with different propagation velocities are mixed together. However, the case $l_0^+=l_0^-$ is excluded even if we restrict to cavities. If $l_0^+\not=l_0^-$ is known a priori, then either $l_0^+<l_0^-$ or $l_0^+>l_0^-$. Thus, from Theorem \ref{mixed but separated case}, we see that if the observation time $T$ is reasonably large (i.e., longer than $2l_0/\sqrt{\gamma_0}$), the indicator function indicates the shortest distance and whether the cavities that achieve it are positive or negative cavities.  
On the other hand, if $l_0^+=l_0^-$, Theorem \ref{mixed but separated case} gives no information such as the shortest distance. \\
\indent
Hence, our goal in this paper is to extract information like Theorem \ref{mixed but separated case} even in the case $l_0^+=l_0^-$. 
For this purpose we need the exact form of the top terms of the indicator function. Therefore, we will analyze the indicator function by using an asymptotic solution under the \lq\lq non-degenerate condition''(cf.\,Definition \ref{non-degenerate condtion for the minal points}).  
For the construction of the asymptotic solution, we put more assumption for the smoothness on the boundaries of the cavities and on the function $\lambda_1$.  We clarify how the curvature of the boundaries and the value of $\lambda_1$ are related to the top term of the expansion of the indicator function. It turns out that information about the cavities can be obtained even for the case $l_0^+=l_0^-$ as Theorem \ref{mixed but separated case} (See the details of the result in section \ref{Main result}).  
\par
In section \ref{Non-degenerate condition}, we give a system of local coordinates around the minimizing points attaining the shortest length, and introduce the non-degenerate condition.  In section \ref{Main result}, we state the main result of this paper. Moreover, we make a reduction for the indicator function to a more manageable form. 
In section \ref{Approximate solutions of w_0}, we construct an asymptotic solution in order to obtain the exact form of the top term of asymptotic expansion of $I_\tau$.  In section \ref{Top term structure}, admitting the estimates of the reminder terms, we analyze the indicator function by the asymptotic solution, eventually obtain the result
(Theorem \ref{a result for the case l_0^+ = l_0^-_general}). In section \ref{Estimates of the remainder terms}, we estimates the remainder terms of this approximation.

\setcounter{equation}{0}
\section{Non-degenerate condition}
\label{Non-degenerate condition}

In addition to (\ref{shortest lenghs1}), another shortest length is defined as follows:
\begin{align*}
l_1 &= \min_{x \in \overline{D}, y, \tilde{y} \in \overline{B}}L(x, y, \tilde{y}), 
\text{ where $L(x, y, \tilde{y})=\v{x - y}+\v{x - \tilde{y}}$. }
\end{align*}
The minimizing points attaining 
$\inf_{(x, y, \tilde{y}) \in \partial{D}{\times}B{\times}B}L(x, y, \tilde{y})$ contribute to the leading part of $I_{\tau}$.  In this article, we call the points attaining $l_0$ or $l_1$ the stationary points. 
As in the proof of Lemma 5.1 in \cite{transmission No1} 
and Lemma 2.5 in \cite{M. and W. Kawashita no1}, 
we see the following facts about the stationary points.
\begin{Lemma}\label{properties of the set of the minimizers}
We have the following properties:
\par\noindent
(1) If $x_0 \in \overline{D}$ and $y_0 \in \overline{B}$ satisfy 
$L_0(x_0,y_0) = l_0$, 
then 
$L(x_0, y_0, y_0) = 2l_0$.  
Hence, we have $l_1 = 2l_0$.  
\par\noindent
(2) If $x_0 \in \overline{D}$ and $y_0 \in \overline{B}$ satisfy 
$L_0(x_0,y_0) = l_0$, 
then $x_0 \in \partial{D}$, $y_0 \in \partial{B}$.
Further, if $\partial{D}$ (resp. $\partial{B}$) has $C^1$ boundary, 
$$
\nu_{x_0} = \frac{y_0 - x_0}{\v{y_0 - x_0}}, 
\qquad \left( \text{resp. } \quad
\nu_{y_0} = -\frac{y_0 - x_0}{\v{y_0 - x_0}} \, \right).
$$
\par\noindent
(3) Suppose that ${D}$ has $C^1$ boundary or $B$ is convex. Then, 
if $x_1 \in \overline{D}$, $y_1, \tilde{y}_1 \in \overline{B}$ satisfy $L(x_1, y_1, \tilde{y}_1) = l_1$,
we have $y_1 = \tilde{y}_1$ and $L_0(x_1,y_1) = l_0$.  
\end{Lemma}

Put $E = \{(x_0, y_0, \tilde{y}_0) \in \overline{D}\times\overline{B}\times\overline{B} \mid 
L(x_0, y_0, \tilde{y}_0) = 2l_0\, \}$, 
$E_0 = \{ (x_0, y_0) \in \overline{D}\times\overline{B} \mid
\v{x_0 - y_0} = l_0\, \}$ and 
$\Gamma=\{ x_0\in \overline{D} \mid
\v{x_0 - y_0} = l_0 \text{ for some }y_0 \in \overline{B}\, \}$. 
From Lemma \ref{properties of the set of the minimizers}, we have 
\begin{equation}
E = \{(x_0, y_0, y_0) \in \partial{D}\times\partial{B}\times\partial{B} \mid (x_0, y_0) \in E_0\}.
\label{minimalpoint}
\end{equation}
The sets $E$, $E_0$ and $\Gamma$ are divided into the following 
three parts: $E = \cup_{\alpha \in\{n_+, n_-, d\}}E^{\alpha}$,  
$E_0 = \cup_{\alpha \in \{n_+, n_-, d\}}E_0^{\alpha}$ and $\Gamma= \cup_{\alpha \in \{ n_+, n_-, d\}}\Gamma^\alpha$, where 
\begin{align*}
&E^{\alpha}= \{ (x_0, y_0, {y}_0) \in E \mid x_0 \in \partial{D^\alpha} \}, 
\quad 
E_0^{\alpha} = \{ (x_0, y_0) \in E_0 \mid x_0 \in \partial{D^\alpha} \}, \\
&\Gamma^{\alpha}= \{ x_0 \in \partial{D^{\alpha}} \mid (x_0, y_0) \in E_0^{\alpha} \text{ for some } y_0 \in \partial{B}\,\}. %
\end{align*}
\par
Let $B_r(a) = \{ x \in \R^3 \vert \v{x - a} < r\}$ for $r > 0$ and $a \in \R^3$.
In this paper, we discuss the problem under the following non-degenerate condition: 
\begin{Def}\label{non-degenerate condtion for the minal points}(non-degenerate condition) \\
We say that $B$ and $D$ satisfy non-degenerate condition for $L_0(x,y)$, if every $(x_0, y_0) \in E_0$ is non-degenerate critical point for $L_0(x,y)$; there exist constants $c_0 > 0$ and $\delta > 0$ satisfying 
$L_0(x, y) \geq l_0 + c_0(\v{x - x_0}^2+\v{y-y_0}^2)$ for 
$(x, y) \in (\partial{D}{\cap}B_\delta(x_0))\times(\partial{B}{\cap}B_\delta(y_0))$. 
\end{Def}

If $B$ and $D$ satisfy the non-degenerate condition for $L_0(x,y)$, $E_0^{\alpha}$ consists only of a finite number ($=M^{\alpha}$) of isolated points: 
$E_0^{\alpha} = \{(x_{j}^{\alpha}, y_{j}^{\alpha}) \in E^{\alpha} \mid j = 1, 2, \ldots, M^{\alpha}\}$ ($\alpha \in \{n_+, n_-, d\}, E_0^\alpha \neq \emptyset$).  
Put $M^{n_+}+M^{n_-}+M^{d} = M (\geq 1)$. If $E_0^\alpha = \emptyset$, 
then we put $M^\alpha = 0$.  In what follows, to express a point in $E_0^\alpha$, we use the notation $(x_0^\alpha, y_0^\alpha) \in E_0^\alpha$ instead of $(x_j^\alpha, y_j^\alpha)$ if we do not need to discriminate them.
\par
%
%
%
%
Hereafter, we assume that both $\partial{D}$ and $\partial{B}$ are $C^2$ class. 
For any point $(x_0^\alpha, y_0^\alpha) \in E_0^\alpha$, we take the following system of local coordinates around $(x_0^\alpha, y_0^\alpha) \in E_0^\alpha$: 
for open neighborhoods $U_{x_0^\alpha}$ and $U_{y_0^\alpha}$ of $(0,0)$ in $\mathbb R^2$ 
there exist $r_0>0$, functions $g^{\alpha}=g_{x_0^\alpha}^\alpha$ and $h^\alpha=h^\alpha_{y_0^\alpha} 
\in C^{2}(\mathbb R^2)$ such that 
$g^{\alpha}(0,0)=h^\alpha(0,0)=0$, 
$\nabla g^{\alpha}(0,0)=\nabla h^\alpha(0,0)=0$ 
and the maps
\begin{align*}
U_{x_0^\alpha}\ni\,\sigma=(\sigma_1,\sigma_2)\mapsto
s^{\alpha}(\sigma)
&=x_0^\alpha+\sigma_1 e_1+\sigma_2 e_2-g^{\alpha}(\sigma_1,\sigma_2)\nu_{x_0^\alpha} 
\\
&\in \partial D^{\alpha}\cap B_{4r_0}(x_0^\alpha), \\
U_{y_0^\alpha}\ni\,u=(u_1, u_2)\mapsto
b^\alpha(u)&=x_0^\alpha+u_1 e_1+u_2 e_2+\{ l_0 + h^\alpha(u_1,u_2)\}\nu_{x_0^\alpha} \\
&\in \partial B\cap B_{4r_0}(y_0^\alpha) 
\end{align*}
give a system of local coordinates around $(x_0^\alpha, y_0^{\alpha})$, where 
$\{e_1, e_2\}$ is an orthogonal basis for $T_{x_0^\alpha}(\partial D)$ and 
$T_{y_0^\alpha}(\partial B)$, $s^{\alpha}(0)=x_0^\alpha$ and $b^\alpha(0)=y_0^\alpha$. 
We can take $r_0 > 0$ small enough so that 
$B_{4r_0}(x_j^\alpha)\cup B_{4r_0}(y_j^\alpha)$ $\alpha \in \{n_+, n_-, d\}$, $j = 1, 2, \ldots, M^{\alpha}$ are mutually disjoint sets. Hence, for any point $(x_0^\alpha, y_0^\alpha) \in E_0^\alpha$, we have $(B_{4r_0}(x_0^\alpha){\times}B_{4r_0}(y_0^\alpha)){\cap}E_0 = \{(x_0^\alpha, y_0^\alpha)\}$. In what follows, we can assume that $\min\{{\rm dist}(D^\alpha, D^\beta) \mid \alpha, \beta \in \{n_+, n_-, d\}, \alpha \neq \beta \} > 9r_0$ and $\min\{{\rm dist}(D^\alpha, B) \mid \alpha \in \{n_+, n_-, d\}\} > 9r_0$ since we can take $r_0 > 0$ smaller if necessary.\\
\indent
For each $\alpha \in \{n_+, n_-, d\}$, we set ${\mathcal U}^\alpha_k = \cup_{j = 1}^{M^{\alpha}}B_{kr_0}(x_j^\alpha)$ ($k = 1, 2, 3$). We also set  
${\mathcal U}_k = \cup_{\alpha \in \{n_+, n_-, d\}}{\mathcal U}^\alpha_k$ ($k = 1, 2, 3$), 
and ${\mathcal V}_1=\cup_{\alpha \in \{n_+, n_-, d\}} \cup_{j = 1}^{M^{\alpha}}B_{r_0}(y_j^\alpha)$.  
Note that $\overline{{\mathcal U}^\alpha_1} \subset {\mathcal U}^\alpha_2$, $\overline{{\mathcal U}^\alpha_2} \subset {\mathcal U}^\alpha_3$, $\overline{{\mathcal U}_1} \subset {\mathcal U}_2$ and $\overline{{\mathcal U}_2} \subset {\mathcal U}_3$, and there exists a constant $c_0 > 0$ such that 
\begin{equation}
\v{x - y} \geq l_0+c_0 
\enskip \text{ if $(x, y) \in (\partial{D}\setminus{\mathcal U}_1)\times\overline{B}$\,  or  \,$(x, y) \in \partial{D} \times(\overline{B}\setminus{\mathcal V}_1)$, } 
\label{estimates on the outside of E^alpha}
\end{equation}
since $E_0 \subset  {\mathcal U}_1\times{\mathcal V}_1$. 
\par
From (\ref{minimalpoint}), the stationary points of $L(x, y, \tilde{y})$ is given by $(x_0^\alpha, y_0^\alpha, y_0^\alpha) \in E^\alpha$ for $(x_0^\alpha, y_0^\alpha) \in E_0^\alpha$.  Using the above local coordinates around $(x_0^\alpha, y_0^{\alpha})$, we define the following notation: 
$$
\tilde{L}^{\alpha}(\sigma, u, \tilde{u}) = L(s^{\alpha}(\sigma), b^{\alpha}(u), b^{\alpha}(\tilde{u})) 
= \v{s^{\alpha}(\sigma)-b^{\alpha}(u)}
+\v{s^{\alpha}(\sigma)-b^{\alpha}(\tilde{u})}. 
$$
\begin{Lemma}\label{equivalent of non-degenerate}
If $B$ is convex, then 
there is an equivalence between 
$(x_0^\alpha, y_0^\alpha)$ being a non-degenerate critical point for $L_0(x,y)$ (Definition \ref{non-degenerate condtion for the minal points}) and $(0, 0, 0)$ being a non-degenerate critical point for $\tilde{L}^{\alpha}(\sigma, u, \tilde{u})$. 
\end{Lemma}
{\it Proof.} Since both $\partial{D}$ and  $\partial{B}$ are $C^2$ class, 
$(x_0^\alpha, y_0^\alpha)$ being a non-degenerate critical point for $L_0(x,y)$ is equivalent to $\text{Hess}(\tilde{L}^\alpha_0)(0,0) > 0$. Similarly, 
$\text{Hess}(\tilde{L}^{\alpha})(0,0,0) > 0$ is equivalent to $\tilde{L}^{\alpha}(\sigma, u, \tilde{u})$ being non-degenerate at $(0, 0, 0)$.  
Then, we calculate $\text{Hess}(\tilde{L}^{\alpha})(0,0,0)$ to prove Lemma \ref{equivalent of non-degenerate}.  
If we set $k^{\alpha}(u, \sigma)=g^\alpha(\sigma)+h^\alpha(u)$, 
then we have $s^\alpha(\sigma)-b^\alpha(u)
=(\sigma_1-u_1)e_1+(\sigma_2-u_2)e_2-(l_0+k^\alpha(u,\sigma))\nu_{x_0^\alpha}$ and 
\begin{align*}
\frac{\partial^2}{\partial \sigma_j \partial \sigma_i} \tilde{L}^\alpha
=& \frac{\delta_{ij}+\frac{\partial^2k^\alpha(u,\sigma)}{\partial \sigma_j \partial \sigma_i}
(l_0 + k^\alpha(u,\sigma))+\frac{\partial k^\alpha(u,\sigma)}{\partial \sigma_i}
\frac{\partial k^\alpha(u,\sigma)}{\partial \sigma_j}}{\v{s^\alpha(\sigma)-b^\alpha(u)}} \\
&-\frac{(\sigma_i-u_i+\frac{\partial k^\alpha(u,\sigma)}{\partial \sigma_i}(l_0 + k^\alpha(u,\sigma)))
(\sigma_j-u_j+\frac{\partial k^\alpha(u,\sigma)}{\partial \sigma_j}(l_0 +k^\alpha(u,\sigma)))}
{\v{s^\alpha(\sigma)-b^\alpha(u)}^3}\\
&+\frac{\delta_{ij}+\frac{\partial^2k^\alpha(\tilde{u},\sigma)}{\partial \sigma_j \partial \sigma_i}
(l_0 +k^\alpha(\tilde{u},\sigma))+\frac{\partial k^\alpha(\tilde{u},\sigma)}{\partial \sigma_i}
\frac{\partial k^\alpha(\tilde{u},\sigma)}{\partial \sigma_j}}{\v{s^\alpha(\sigma)-b^\alpha(\tilde{u})}} \\
&-\frac{(\sigma_i-\tilde{u}_i+\frac{\partial k^\alpha(\tilde{u},\sigma)}{\partial \sigma_i}(l_0 +k^\alpha(\tilde{u},\sigma)))
(\sigma_j-\tilde{u}_j+\frac{\partial k^\alpha(\tilde{u},\sigma)}{\partial \sigma_j}(l_0 +k^\alpha(\tilde{u},\sigma)))}
{\v{s^\alpha(\sigma)-b^\alpha(\tilde{u})}^3},
\end{align*}
\begin{align*}
\frac{\partial^2}{\partial \sigma_j \partial u_i} \tilde{L}^\alpha
=& \frac{-\delta_{ij}+\frac{\partial k^\alpha(u,\sigma)}{\partial u_i}
\frac{\partial k^\alpha(u,\sigma)}{\partial \sigma_j}}{\v{s^\alpha(\sigma)-b^\alpha(u)}} \\
&-\frac{(u_i-\sigma_i+\frac{\partial k^\alpha(u,\sigma)}{\partial u_i}(l_0 +k^\alpha(u,\sigma)))
(\sigma_j-u_j+\frac{\partial k^\alpha(u,\sigma)}{\partial \sigma_j}
(l_0 + k^\alpha(u,\sigma)))}
{\v{s^\alpha(\sigma)-b^\alpha(u)}^3}, \\
\frac{\partial^2}{\partial \sigma_j \partial \tilde{u}_i} \tilde{L}^\alpha
=& \frac{-\delta_{ij}+\frac{\partial k^\alpha(\tilde{u},\sigma)}{\partial \tilde{u}_i}
\frac{\partial k^\alpha(\tilde{u},\sigma)}{\partial \sigma_j}}{\v{s^\alpha(\sigma)-b^\alpha(\tilde{u})}} \\
&-\frac{(\tilde{u}_i-\sigma_i+\frac{\partial k^\alpha(\tilde{u},\sigma)}{\partial \tilde{u}_i}
(l_0 +k^\alpha(\tilde{u},\sigma)))
(\sigma_j-\tilde{u}_j+\frac{\partial k^\alpha(\tilde{u},\sigma)}{\partial \sigma_j}
(l_0 +k^\alpha(\tilde{u},\sigma)))}
{\v{s^\alpha(\sigma)-b^\alpha(\tilde{u})}^3}. 
\end{align*}
Since $\tilde{L}^\alpha(\sigma, u, \tilde{u})$ 
takes a extreme value $2l_0$ at $(0, 0, 0)$, we know that 
$\frac{\partial k^\alpha(0,0)}{\partial \sigma_j}
=\frac{\partial k^\alpha(0,0)}{\partial {u}_j}
=0$. 
Set the $2 \times 2$ matrices $G^{\alpha}=(g_{ij}^{\alpha})$, $H^\alpha=(h_{ij}^\alpha)$, 
$S_{g^{\alpha}}=I+l_0 G^{\alpha}$ and $S_{h^\alpha}=I+l_0 H^\alpha$,  
where $g_{ij}^\alpha=\frac{\partial^2 g^{\alpha}(0)}{\partial \sigma_j\partial\sigma_i}$, 
$h_{ij}^\alpha=\frac{\partial^2 h^\alpha(0) }{\partial u_j\partial u_i}$. 
Then, at the point $(\sigma, u, \tilde{u})=(0,0,0)$, 
the $6 \times 6$ matrix $\text{Hess}(\tilde{L}^\alpha)(\sigma, u, \tilde{u})$ has the form  
\begin{align*}
\text{Hess}(\tilde{L}^\alpha)(0,0,0)=\frac{1}{l_0}
\left(
\begin{array}{ccc}
2S_{g^{\alpha}} & -I & -I \\
-I & S_{h^\alpha} & O \\
-I & O & S_{h^\alpha} 
\end{array}
\right).
\end{align*}
Set $\tilde{L}^\alpha_0(\sigma, u) = L_0(s^\alpha(\sigma), b^\alpha(u)) = \v{s^\alpha(\sigma) - b^\alpha(u)}$, then we have 
$$
\text{Hess}(\tilde{L}^\alpha_0)(0,0)=\frac{1}{l_0}
\left(
\begin{array}{cc}
S_{g^{\alpha}} & -I \\
-I & S_{h^\alpha}  
\end{array}
\right)
$$
in the same way as above.  Thus, if we set 
$P = \begin{pmatrix}
I & O & O \\
O & 2^{-1}I & 2^{-1}I \\
O & 2^{-1}I & -2^{-1}I
\end{pmatrix}$
we have
\begin{align}
\text{Hess}(\tilde{L}^\alpha)(0,0,0)=\frac{2}{l_0}{}^t{\hskip-0.5mm}P\begin{pmatrix}
l_0\text{Hess}(\tilde{L}^\alpha_0)(0,0) & O \\ O & S_{h^\alpha}
\end{pmatrix}
P. 
\label{decomposition_of_matrix}
\end{align}
If $B$ is convex, we have $S_{h^\alpha} > 0$ because $H^\alpha \geq 0$.  
Therefore, when $B$ is convex, 
$\text{Hess}(\tilde{L}^\alpha_0)(0,0) > 0$ and 
$\text{Hess}(\tilde{L}^{\alpha})(0,0,0) > 0$ are equivalent, which implies Lemma \ref{equivalent of non-degenerate}. \\
\hfill$\square$ 
\par
To obtain the asymptotic behavior of $I_{\tau}$, the exact form of 
$\text{det}( \text{Hess}(\tilde{L}^\alpha)(0,0,0))$ is needed. 
We denote the mean and Gauss curvatures of $\partial{B}$ at $y \in \partial{B}$ by ${\mathcal M}_{\partial B}(y)$ and ${\mathcal G}_{\partial B}(y)$, respectively. Then, from (\ref{decomposition_of_matrix}) it follows that 
\[
\text{det}(\text{Hess}(\tilde{L}^\alpha)(0,0,0)) = \frac{4}{(l_0)^4}{\mathcal A}_\alpha(x_0^\alpha; \partial{D}, \partial{B}), 
\]
where
\begin{align*}
{\mathcal A}_\alpha(x_0^\alpha; \partial{D}, \partial{B}) = \big(1&+2l_0{\mathcal M}_{\partial B}(y_0^\alpha)
+(l_0)^2{\mathcal G}_{\partial B}(y_0^\alpha) \big)
\big\{(l_0)^2{\mathcal G}_{\partial D}(x_0^\alpha){\mathcal G}_{\partial B}(y_0^\alpha) 
\\
&+ 2l_0\big({\mathcal M}_{\partial B}(y_0^\alpha){\mathcal G}_{\partial D}(x_0^\alpha) + {\mathcal M}_{\partial D}(x_0^\alpha){\mathcal G}_{\partial B}(y_0^\alpha)\big) + \det(G^\alpha + H^\alpha)\big\}.
\end{align*}
For simplicity, we take $B$ as a ball with radius $a > 0$, then we have
\begin{align*}
&{\mathcal M}_{\partial{B}}(y_0^\alpha) = \frac{1}{a}, \enskip\enskip {\mathcal G}_{\partial{B}}(y_0^\alpha) = \frac{1}{a^2}, \enskip\enskip
H^\alpha = \frac{1}{a}I, \\ 
&\det(G^\alpha+H^\alpha) 
= \frac{1}{a^2}\det(I+aG^\alpha) =  \frac{1}{a^2}\Big(1+2a{\mathcal M}_{\partial{D}}(x_0^\alpha) + a^2{\mathcal G}_{\partial{D}}(x_0^\alpha)\Big). 
\end{align*}
Thus, 
\begin{align*}
{\mathcal A}_\alpha(x_0^\alpha; \partial{D}, \partial{B}) = \left(\frac{a+l_0}{a}\right)^4\left({\mathcal G}_{\partial D}(x_0^\alpha) + 2\frac{{\mathcal M}_{\partial{D}}(x_0^\alpha)}{a+l_0} + \frac{1}{(a+l_0)^2}\right).
\end{align*}
%
Let us denote the principal curvatures of $\partial{D^{\alpha}}$ at $x \in \partial{D^{\alpha}}$ 
by $\kappa_{1, \alpha}(x)$ and $\kappa_{2, \alpha}(x)$ 
with $\kappa_{1, \alpha}(x) \leq \kappa_{2, \alpha}(x)$ for $\alpha \in \{n_+, n_-, d\}$. Then, it follows that $2{\mathcal M}_{\partial{D}}(x_0^\alpha)=\kappa_{1, \alpha}(x_0^\alpha)+\kappa_{2, \alpha}(x_0^\alpha)$ and 
${\mathcal G}_{\partial D}(x_0^\alpha)=\kappa_{1, \alpha}(x_0^\alpha)\kappa_{2, \alpha}(x_0^\alpha)$.  
Hence, we have
\begin{align*}
\text{det}( \text{Hess}&(\tilde{L}^\alpha)(0,0,0)) 
= \frac{4(a+l_0)^4}{(l_0)^4a^4}{\mathcal A}_\alpha(x_0^\alpha),
\end{align*}
where
${\mathcal A}_\alpha(x) = \prod_{j = 1}^2\Big(\kappa_{j, \alpha}(x) + \frac{1}{l_0 + a}\Big)$ ($\alpha \in \{n_+, n_-, d\}$). 
From the non-degenerate condition (Definition \ref{non-degenerate condtion for the minal points}) it follows that 
\begin{align*}
{\mathcal A}_{\alpha}(x) > 0 \quad
\text{ for $(x, y) \in E_0^{\alpha}, \alpha \in\{n_+, n_-, d\}$. }
\end{align*}
\par
\setcounter{equation}{0}
\section{Main results and a reduction of $I_{\tau}$} 
\label{Main result}
The main result of this paper is as follows:
\begin{Thm}
\label{a result for the case l_0^+ = l_0^-_general}
Assume that $\partial{D^n}$ and $\partial{D^d}$ are $C^4$ class, ${B}$ is convex set with $C^2$ boundary, 
$\lambda_0 \in L^{\infty}(\partial D^n)$, $\lambda_1 \in C^2(\partial D^n)$ and $f \in C^1(\overline{B})$, moreover, $D$ and $B$ satisfy non-degenerate condition for $L_0(x,y)$. Then, 
there exists $\delta > 0$ such that 
\begin{align*}
I_{\tau} = {\frac{\pi\gamma_0}{\tau^4}}\hskip-2pt\sum_{\alpha \in \{ n_+, n_-, d \}}
\sum_{j=1}^{M^{\alpha}}e^{-\frac{{\tau}}{\sqrt{\gamma_0}}2l_0}&\Bigg\{b_\alpha({x_j^{\alpha}})
\frac{(f({y_j^{\alpha}}))^2}{2\sqrt{{\mathcal A}_\alpha({x_j^{\alpha}}; \partial{D}, \partial{B})}}
+ O(\tau^{-\frac{1}{2}})\Bigg\}
\\&
+ O(e^{-\frac{\tau}{\sqrt{\gamma_0}}(2l_0+\delta)}) + O(\tau^{-1}e^{-\tau{T}})
\quad(\tau \to \infty),
\end{align*}
where $(x_j^{\alpha}, y_j^{\alpha}) \in E_0^{\alpha}$, 
$b_{\alpha}(x) 
= \frac{\sqrt{\gamma_0}-\lambda_1(x)}
{\sqrt{\gamma_0}+\lambda_1(x)}$ for $\alpha\in \{n_+, n_-\}$ and 
$b_{d}(x) = -1$. \par
If we set ${\mathcal{T}_0} =\sum_{\alpha \in \{ n_+, n_-, d \}}
\sum_{j=1}^{M^{\alpha}} b_\alpha({x_j^{\alpha}})
\frac{(f({y_j^{\alpha}}))^2}{2\sqrt{{\mathcal A}_\alpha({x_j^{\alpha}}; \partial{D}, \partial{B})}}$, then   
\par\noindent
(1) for $T < 2l_0/\sqrt{\gamma_0}$, $\lim_{\tau \to \infty}e^{{\tau}T}I_\tau = 0$, 
\par\noindent
(2) for $T > 2l_0/\sqrt{\gamma_0}$, 
if (\ref{emission condition}) holds, then we have: 
$$\text{$(ii)_{+0}$} \quad
\lim_{\tau \to \infty}e^{{\tau}T}I_\tau = +\infty \enskip \text{if ${\mathcal{T}_0} >0$,} 
\qquad 
\text{$(ii)_{-0}$} \quad
\lim_{\tau \to \infty}e^{{\tau}T}I_\tau =-\infty \enskip \text{if ${\mathcal{T}_0}<0$.} 
$$
Further, in each of the cases $(ii)_{+0}$ and $(ii)_{-0}$, we also have
\begin{equation*}
\lim_{\tau \to \infty}\frac{\sqrt{\gamma_0}}{2\tau}\log\v{I_\tau} = 
\begin{cases} -l_0^+ & \text{for the case of $(ii)_{+0}$},
\\
-l_0^- & \text{for the case of $(ii)_{-0}$}.
\end{cases}
\end{equation*}
\end{Thm}
\begin{Remark}\label{results of general case} 
From the asymptotic form of $I_{\tau}$ of Theorem \ref{a result for the case l_0^+ = l_0^-_general}, the same information as in Theorem \ref{mixed but separated case} can be obtained without knowing a priori that $l_0^+ \not=l_0^-$. If ${\mathcal{T}_0}$ is not zero, the shortest distance $l_0$ can be obtained from Theorem \ref{a result for the case l_0^+ = l_0^-_general}, and it can be determined whether $l_0$ is $l_0^+$ or $l_0^-$. 
If $l_0^+ =l_0 < l_0^-$, then ${\mathcal{T}_0} >0$, and 
if $l_0^- =l_0 < l_0^+$, then ${\mathcal{T}_0} <0$. 
\end{Remark}
To make the meaning of ${\mathcal{T}_0} >0$ or ${\mathcal{T}_0} <0$ clearer, in particular, take $B$ to be a ball with radius $a > 0$.
Then, Theorem \ref{a result for the case l_0^+ = l_0^-_general} becomes as follows:
\begin{Thm}
\label{a result for the case l_0^+ = l_0^-}
Assume that $\partial{D^n}$ and  $\partial{D^d}$ are $C^4$ class, $B$ is a ball with radius $a > 0$, 
and $D$ and $B$ satisfy non-degenerate condition for $L_0(x,y)$, 
$\lambda_0 \in L^{\infty}(\partial D^n)$, $\lambda_1 \in C^2(\partial D^n)$ and $f \in C^1(\overline{B})$.  
Then, there exists $\delta > 0$ such that 
\begin{align*}
I_{\tau} = {\frac{\pi\gamma_0}{\tau^4}}
&e^{-\frac{{\tau}}{\sqrt{\gamma_0}}2l_0}
\Big\{ {\mathcal{T}_0} + O(\tau^{-\frac{1}{2}})\Big\}+ O(e^{-\frac{\tau}{\sqrt{\gamma_0}}(2l_0+\delta)}) + O(\tau^{-1}e^{-\tau{T}})
\quad(\tau \to \infty), 
\end{align*}
where ${\mathcal{T}_0}=\sum_{\alpha \in \{ n_+, n_-, d \}}
\sum_{j=1}^{M^{\alpha}}\frac{a^2}{2(l_0 + a)^2}b_\alpha({x_j^{\alpha}})
\frac{(f({y_j^{\alpha}}))^2}{\sqrt{{\mathcal A}_\alpha({x_j^{\alpha}})}}$. 
\end{Thm}
As stated in Theorem \ref{a result for the case l_0^+ = l_0^-_general}, in the case $l_0^+ \not =l_0^-$, the asymptotic form of $I_{\tau}$ obtained in Theorem \ref{a result for the case l_0^+ = l_0^-} provides the same information as in Theorem \ref{mixed but separated case}. Furthermore, for the case $l_0^+ =l_0^-$, which is excluded in Theorem \ref{mixed but separated case}, information on the shortest distance can also be extracted. Let us consider two simple examples of the case $l_0^+ =l_0^-$ in order to explain Theorem \ref{a result for the case l_0^+ = l_0^-}.

\begin{rei}\label{lambda1=0}
Consider the case that $\lambda_1(x) = 0$ $(x \in \partial{D}^n)$ and only two pairs of the points $(x_0^{n_+}, y_0^{n_+})$ and $(x_0^{d}, y_0^{d})$ attain $l_0$.  
In this case, 
we have $D^{n_-} = \emptyset$, and $E_0^{n_+}$ and  $E_0^{d}$ are expressed by 
$E_0^{n_+} = \{(x_0^{n_+}, y_0^{n_+})\}$ and $E_0^{d} = \{(x_0^{d}, y_0^{d})\}$.
If $f(y) = 1 $ in $\overline{B}$, Theorem \ref{a result for the case l_0^+ = l_0^-} implies that 
$$
\mathcal{T}_0=\frac{a^2}{2(l_0+a)^2}\Big(\frac{1}{\sqrt{{\mathcal A}_{n_+}(x_0^{n_+})}} - \frac{1}{\sqrt{{\mathcal A}_{d}(x_0^{d})}}\Big).
$$
Thus, the indicator function $I_\tau$ detects the boundary which the relative Gaussian curvature is smaller. 
Considering the case in Figure 1, $\lim_{\tau \to \infty} e^{\tau T}I_{\tau}=\infty$ when $T >2l_0$, the indicator function $I_\tau$ detects $D^{n_+}$.
 \begin{figure}[ht]
 \hskip20mm
{\unitlength 0.1in%
\begin{picture}(37.2300,18.4100)(18.9000,-24.5100)%
%
\special{pn 8}%
\special{ar 4118 905 0 0 0.1023542 0.1218872}%
%
\special{pn 8}%
\special{ar 4205 813 203 203 2.2794226 2.2348425}%
%
\special{pn 8}%
\special{pa 4653 1651}%
\special{pa 4295 981}%
\special{fp}%
\put(41.1000,-8.8000){\makebox(0,0)[lb]{$B$}}%
\put(22.8000,-16.2000){\makebox(0,0)[lb]{$D^d$}}%
\put(46.2000,-20.4000){\makebox(0,0)[lb]{$D^{n_+}$}}%
\put(34.7000,-10.9000){\makebox(0,0)[lb]{$l_0$}}%
\put(45.0000,-12.8000){\makebox(0,0)[lb]{$l_0$}}%
%
\special{pn 8}%
\special{pa 5600 1560}%
\special{pa 5609 1591}%
\special{pa 5613 1622}%
\special{pa 5610 1654}%
\special{pa 5601 1685}%
\special{pa 5590 1715}%
\special{pa 5576 1744}%
\special{pa 5560 1772}%
\special{pa 5543 1799}%
\special{pa 5503 1849}%
\special{pa 5481 1873}%
\special{pa 5459 1896}%
\special{pa 5437 1918}%
\special{pa 5413 1941}%
\special{pa 5390 1962}%
\special{pa 5366 1983}%
\special{pa 5341 2004}%
\special{pa 5316 2023}%
\special{pa 5290 2043}%
\special{pa 5264 2061}%
\special{pa 5238 2080}%
\special{pa 5211 2097}%
\special{pa 5185 2115}%
\special{pa 5158 2132}%
\special{pa 5130 2149}%
\special{pa 5103 2166}%
\special{pa 5075 2182}%
\special{pa 5048 2197}%
\special{pa 5019 2213}%
\special{pa 4991 2227}%
\special{pa 4933 2255}%
\special{pa 4905 2269}%
\special{pa 4875 2283}%
\special{pa 4846 2296}%
\special{pa 4817 2308}%
\special{pa 4787 2321}%
\special{pa 4697 2354}%
\special{pa 4637 2374}%
\special{pa 4606 2384}%
\special{pa 4575 2393}%
\special{pa 4545 2402}%
\special{pa 4483 2418}%
\special{pa 4451 2425}%
\special{pa 4420 2431}%
\special{pa 4388 2436}%
\special{pa 4357 2441}%
\special{pa 4325 2445}%
\special{pa 4293 2448}%
\special{pa 4261 2450}%
\special{pa 4229 2451}%
\special{pa 4197 2451}%
\special{pa 4165 2448}%
\special{pa 4134 2444}%
\special{pa 4102 2439}%
\special{pa 4071 2431}%
\special{pa 4041 2420}%
\special{pa 4013 2404}%
\special{pa 3987 2385}%
\special{pa 3967 2360}%
\special{pa 3952 2331}%
\special{pa 3945 2300}%
\special{pa 3946 2268}%
\special{pa 3953 2236}%
\special{pa 3962 2206}%
\special{pa 3975 2176}%
\special{pa 3990 2148}%
\special{pa 4007 2121}%
\special{pa 4025 2095}%
\special{pa 4045 2070}%
\special{pa 4066 2046}%
\special{pa 4088 2022}%
\special{pa 4156 1954}%
\special{pa 4180 1933}%
\special{pa 4205 1912}%
\special{pa 4229 1892}%
\special{pa 4307 1835}%
\special{pa 4333 1817}%
\special{pa 4360 1799}%
\special{pa 4386 1782}%
\special{pa 4414 1765}%
\special{pa 4441 1748}%
\special{pa 4469 1732}%
\special{pa 4496 1716}%
\special{pa 4525 1701}%
\special{pa 4581 1671}%
\special{pa 4668 1629}%
\special{pa 4726 1603}%
\special{pa 4756 1591}%
\special{pa 4785 1579}%
\special{pa 4815 1567}%
\special{pa 4845 1556}%
\special{pa 4875 1546}%
\special{pa 4906 1536}%
\special{pa 4936 1526}%
\special{pa 4967 1516}%
\special{pa 4998 1507}%
\special{pa 5028 1499}%
\special{pa 5059 1491}%
\special{pa 5091 1484}%
\special{pa 5122 1478}%
\special{pa 5154 1472}%
\special{pa 5185 1467}%
\special{pa 5217 1463}%
\special{pa 5281 1457}%
\special{pa 5313 1455}%
\special{pa 5345 1455}%
\special{pa 5377 1456}%
\special{pa 5409 1459}%
\special{pa 5440 1463}%
\special{pa 5471 1470}%
\special{pa 5502 1480}%
\special{pa 5532 1493}%
\special{pa 5559 1510}%
\special{pa 5581 1533}%
\special{pa 5599 1559}%
\special{pa 5600 1560}%
\special{fp}%
%
\special{pn 4}%
\special{pa 5190 1470}%
\special{pa 4210 2450}%
\special{fp}%
\special{pa 4960 1520}%
\special{pa 4060 2420}%
\special{fp}%
\special{pa 4650 1650}%
\special{pa 3960 2340}%
\special{fp}%
\special{pa 5380 1460}%
\special{pa 4410 2430}%
\special{fp}%
\special{pa 5530 1490}%
\special{pa 4650 2370}%
\special{fp}%
\special{pa 5610 1590}%
\special{pa 4960 2240}%
\special{fp}%
%
\special{pn 4}%
\special{sh 1}%
\special{ar 4650 1630 8 8 0 6.2831853}%
\special{sh 1}%
\special{ar 4660 1640 8 8 0 6.2831853}%
\special{sh 1}%
\special{ar 4660 1640 8 8 0 6.2831853}%
\put(33.9000,-14.1000){\makebox(0,0)[lb]{$x_0^d$}}%
\put(46.3000,-15.7000){\makebox(0,0)[lb]{$x_0^{n_+}$}}%
%
\special{pn 4}%
\special{sh 1}%
\special{ar 4030 920 8 8 0 6.2831853}%
\special{sh 1}%
\special{ar 4030 920 8 8 0 6.2831853}%
\put(37.6000,-8.8000){\makebox(0,0)[lb]{$y_0^d$}}%
\put(44.9000,-9.5000){\makebox(0,0){$y_0^{n_+}$}}%
%
\special{pn 8}%
\special{pa 3156 1587}%
\special{pa 3130 1606}%
\special{pa 3104 1624}%
\special{pa 3078 1643}%
\special{pa 3024 1677}%
\special{pa 2968 1709}%
\special{pa 2912 1739}%
\special{pa 2883 1753}%
\special{pa 2854 1766}%
\special{pa 2825 1780}%
\special{pa 2795 1792}%
\special{pa 2766 1804}%
\special{pa 2736 1815}%
\special{pa 2705 1826}%
\special{pa 2645 1846}%
\special{pa 2614 1856}%
\special{pa 2584 1865}%
\special{pa 2522 1883}%
\special{pa 2491 1890}%
\special{pa 2460 1896}%
\special{pa 2396 1908}%
\special{pa 2365 1913}%
\special{pa 2301 1921}%
\special{pa 2270 1924}%
\special{pa 2238 1926}%
\special{pa 2206 1927}%
\special{pa 2174 1927}%
\special{pa 2142 1925}%
\special{pa 2110 1922}%
\special{pa 2078 1917}%
\special{pa 2047 1911}%
\special{pa 2016 1902}%
\special{pa 1986 1890}%
\special{pa 1958 1875}%
\special{pa 1932 1856}%
\special{pa 1911 1832}%
\special{pa 1896 1803}%
\special{pa 1891 1772}%
\special{pa 1892 1740}%
\special{pa 1900 1709}%
\special{pa 1913 1679}%
\special{pa 1928 1651}%
\special{pa 1946 1624}%
\special{pa 1966 1599}%
\special{pa 1987 1576}%
\special{pa 2010 1553}%
\special{pa 2033 1531}%
\special{pa 2057 1510}%
\special{pa 2082 1489}%
\special{pa 2107 1470}%
\special{pa 2133 1451}%
\special{pa 2160 1434}%
\special{pa 2187 1416}%
\special{pa 2215 1400}%
\special{pa 2242 1384}%
\special{pa 2270 1368}%
\special{pa 2298 1353}%
\special{pa 2327 1338}%
\special{pa 2355 1324}%
\special{pa 2384 1310}%
\special{pa 2414 1298}%
\special{pa 2443 1285}%
\special{pa 2503 1261}%
\special{pa 2532 1250}%
\special{pa 2562 1239}%
\special{pa 2593 1228}%
\special{pa 2623 1219}%
\special{pa 2685 1201}%
\special{pa 2747 1187}%
\special{pa 2779 1180}%
\special{pa 2810 1174}%
\special{pa 2841 1169}%
\special{pa 2873 1163}%
\special{pa 2905 1159}%
\special{pa 2936 1155}%
\special{pa 2968 1152}%
\special{pa 3000 1150}%
\special{pa 3032 1150}%
\special{pa 3096 1152}%
\special{pa 3160 1158}%
\special{pa 3192 1164}%
\special{pa 3223 1172}%
\special{pa 3252 1184}%
\special{pa 3280 1200}%
\special{pa 3305 1220}%
\special{pa 3327 1244}%
\special{pa 3342 1272}%
\special{pa 3350 1303}%
\special{pa 3348 1335}%
\special{pa 3340 1366}%
\special{pa 3328 1396}%
\special{pa 3312 1424}%
\special{pa 3294 1450}%
\special{pa 3274 1476}%
\special{pa 3253 1500}%
\special{pa 3231 1523}%
\special{pa 3183 1565}%
\special{pa 3161 1582}%
\special{fp}%
%
\special{pn 4}%
\special{pa 2754 1184}%
\special{pa 2421 1903}%
\special{fp}%
\special{pa 2854 1168}%
\special{pa 2522 1887}%
\special{fp}%
\special{pa 2955 1153}%
\special{pa 2635 1846}%
\special{fp}%
\special{pa 3050 1150}%
\special{pa 2741 1817}%
\special{fp}%
\special{pa 3139 1160}%
\special{pa 2860 1763}%
\special{fp}%
\special{pa 3222 1183}%
\special{pa 2985 1696}%
\special{fp}%
\special{pa 3299 1218}%
\special{pa 3115 1616}%
\special{fp}%
\special{pa 3346 1318}%
\special{pa 3275 1472}%
\special{fp}%
\special{pa 2647 1213}%
\special{pa 2321 1919}%
\special{fp}%
\special{pa 2534 1254}%
\special{pa 2226 1922}%
\special{fp}%
\special{pa 2422 1296}%
\special{pa 2131 1925}%
\special{fp}%
\special{pa 2303 1350}%
\special{pa 2048 1902}%
\special{fp}%
\special{pa 2173 1430}%
\special{pa 1965 1879}%
\special{fp}%
\special{pa 2025 1548}%
\special{pa 1900 1818}%
\special{fp}%
%
\special{pn 4}%
\special{sh 1}%
\special{ar 4310 1000 8 8 0 6.2831853}%
\special{sh 1}%
\special{ar 4310 1000 8 8 0 6.2831853}%
%
\special{pn 8}%
\special{pa 4020 920}%
\special{pa 3340 1240}%
\special{fp}%
%
\special{pn 4}%
\special{sh 1}%
\special{ar 3340 1240 8 8 0 6.2831853}%
\end{picture}}%
\caption{Example \ref{lambda1=0}}
\label{ex1}
 \end{figure}
\end{rei}
\begin{rei}\label{lambda1not=0}
Consider the case that $\lambda_1(x) > 0$ $(x \in \partial{D}^n)$,  only two pairs of the points $(x_0^{n_+}, y_0^{n_+})$ and $(x_0^{d}, y_0^{d})$ attain $l_0$, moreover $\kappa_{j, n_+}(x_0^{n_+})=\kappa_{j,d}(x_0^d)$ for $j=1,2$ and $f(y)=1$ in $\overline{B}$(see Figire 2).  
From ${\mathcal A}_{n_+}(x_0^{n_+})
={\mathcal A}_{d}(x_0^{d})=:{\mathcal A}_0$ and $\lambda_1(x_0^{n_+}) > 0$, it follows that 

$$\mathcal{T}_0=\frac{-a^2}{(l_0+a)^2\sqrt{\mathcal{A}_0}}\frac{\lambda_1(x_0^{n_+})}{\sqrt{\gamma_0}+\lambda_1(x_0^{n_+})}<0.$$   
Then, we have   
$\lim_{\tau \to \infty} e^{\tau T}I_{\tau}=-\infty$ for $T >2l_0$, namely the indicator function $I_\tau$ detects $D^d$.  
 \begin{figure}[ht]
 \hskip20mm
{\unitlength 0.1in%
\begin{picture}(30.4000,17.1000)(20.6000,-20.8300)%
%
\special{pn 8}%
\special{ar 3572 1747 0 0 2.9593245 2.9793541}%
%
\special{pn 8}%
\special{pa 4498 1096}%
\special{pa 4495 1089}%
\special{fp}%
%
\special{pn 8}%
\special{ar 3550 1881 202 202 2.4601950 2.4133712}%
%
\special{pn 8}%
\special{pa 5032 1274}%
\special{pa 5004 1289}%
\special{pa 4975 1302}%
\special{pa 4944 1312}%
\special{pa 4912 1317}%
\special{pa 4881 1320}%
\special{pa 4848 1320}%
\special{pa 4817 1318}%
\special{pa 4753 1310}%
\special{pa 4722 1304}%
\special{pa 4690 1297}%
\special{pa 4660 1288}%
\special{pa 4629 1277}%
\special{pa 4599 1266}%
\special{pa 4570 1254}%
\special{pa 4541 1241}%
\special{pa 4483 1213}%
\special{pa 4455 1198}%
\special{pa 4427 1182}%
\special{pa 4399 1165}%
\special{pa 4372 1148}%
\special{pa 4346 1131}%
\special{pa 4294 1093}%
\special{pa 4244 1053}%
\special{pa 4196 1011}%
\special{pa 4173 988}%
\special{pa 4129 942}%
\special{pa 4108 917}%
\special{pa 4088 893}%
\special{pa 4069 867}%
\special{pa 4049 842}%
\special{pa 4031 815}%
\special{pa 4014 789}%
\special{pa 3998 761}%
\special{pa 3984 732}%
\special{pa 3972 702}%
\special{pa 3963 671}%
\special{pa 3955 640}%
\special{pa 3950 609}%
\special{pa 3948 577}%
\special{pa 3949 545}%
\special{pa 3956 513}%
\special{pa 3968 483}%
\special{pa 3986 457}%
\special{pa 4009 435}%
\special{pa 4034 415}%
\special{pa 4063 400}%
\special{pa 4093 390}%
\special{pa 4125 384}%
\special{pa 4157 381}%
\special{pa 4189 380}%
\special{pa 4221 381}%
\special{pa 4252 383}%
\special{pa 4284 387}%
\special{pa 4316 393}%
\special{pa 4347 401}%
\special{pa 4377 410}%
\special{pa 4408 420}%
\special{pa 4438 431}%
\special{pa 4468 443}%
\special{pa 4498 454}%
\special{pa 4527 467}%
\special{pa 4555 482}%
\special{pa 4584 497}%
\special{pa 4611 513}%
\special{pa 4639 529}%
\special{pa 4666 546}%
\special{pa 4693 564}%
\special{pa 4745 600}%
\special{pa 4771 620}%
\special{pa 4795 640}%
\special{pa 4819 661}%
\special{pa 4843 683}%
\special{pa 4866 705}%
\special{pa 4889 728}%
\special{pa 4933 774}%
\special{pa 4975 822}%
\special{pa 4995 848}%
\special{pa 5012 874}%
\special{pa 5044 930}%
\special{pa 5059 958}%
\special{pa 5073 987}%
\special{pa 5084 1017}%
\special{pa 5092 1048}%
\special{pa 5098 1080}%
\special{pa 5100 1112}%
\special{pa 5099 1144}%
\special{pa 5094 1175}%
\special{pa 5085 1206}%
\special{pa 5069 1234}%
\special{pa 5047 1258}%
\special{pa 5039 1266}%
\special{fp}%
%
\special{pn 8}%
\special{pa 5032 1274}%
\special{pa 5004 1289}%
\special{pa 4975 1302}%
\special{pa 4944 1312}%
\special{pa 4912 1317}%
\special{pa 4881 1320}%
\special{pa 4848 1320}%
\special{pa 4817 1318}%
\special{pa 4753 1310}%
\special{pa 4722 1304}%
\special{pa 4690 1297}%
\special{pa 4660 1288}%
\special{pa 4629 1277}%
\special{pa 4599 1266}%
\special{pa 4570 1254}%
\special{pa 4541 1241}%
\special{pa 4483 1213}%
\special{pa 4455 1198}%
\special{pa 4427 1182}%
\special{pa 4399 1165}%
\special{pa 4372 1148}%
\special{pa 4346 1131}%
\special{pa 4294 1093}%
\special{pa 4244 1053}%
\special{pa 4196 1011}%
\special{pa 4173 988}%
\special{pa 4129 942}%
\special{pa 4108 917}%
\special{pa 4088 893}%
\special{pa 4069 867}%
\special{pa 4049 842}%
\special{pa 4031 815}%
\special{pa 4014 789}%
\special{pa 3998 761}%
\special{pa 3984 732}%
\special{pa 3972 702}%
\special{pa 3963 671}%
\special{pa 3955 640}%
\special{pa 3950 609}%
\special{pa 3948 577}%
\special{pa 3949 545}%
\special{pa 3956 513}%
\special{pa 3968 483}%
\special{pa 3986 457}%
\special{pa 4009 435}%
\special{pa 4034 415}%
\special{pa 4063 400}%
\special{pa 4093 390}%
\special{pa 4125 384}%
\special{pa 4157 381}%
\special{pa 4189 380}%
\special{pa 4221 381}%
\special{pa 4252 383}%
\special{pa 4284 387}%
\special{pa 4316 393}%
\special{pa 4347 401}%
\special{pa 4377 410}%
\special{pa 4408 420}%
\special{pa 4438 431}%
\special{pa 4468 443}%
\special{pa 4498 454}%
\special{pa 4527 467}%
\special{pa 4555 482}%
\special{pa 4584 497}%
\special{pa 4611 513}%
\special{pa 4639 529}%
\special{pa 4666 546}%
\special{pa 4693 564}%
\special{pa 4745 600}%
\special{pa 4771 620}%
\special{pa 4795 640}%
\special{pa 4819 661}%
\special{pa 4843 683}%
\special{pa 4866 705}%
\special{pa 4889 728}%
\special{pa 4933 774}%
\special{pa 4975 822}%
\special{pa 4995 848}%
\special{pa 5012 874}%
\special{pa 5044 930}%
\special{pa 5059 958}%
\special{pa 5073 987}%
\special{pa 5084 1017}%
\special{pa 5092 1048}%
\special{pa 5098 1080}%
\special{pa 5100 1112}%
\special{pa 5099 1144}%
\special{pa 5094 1175}%
\special{pa 5085 1206}%
\special{pa 5069 1234}%
\special{pa 5047 1258}%
\special{pa 5039 1266}%
\special{fp}%
%
\special{pn 8}%
\special{ar 2420 800 790 390 6.1822792 6.1835167}%
%
\special{pn 8}%
\special{pa 3319 497}%
\special{pa 3327 528}%
\special{pa 3330 560}%
\special{pa 3325 591}%
\special{pa 3317 623}%
\special{pa 3306 653}%
\special{pa 3292 682}%
\special{pa 3276 709}%
\special{pa 3258 736}%
\special{pa 3218 786}%
\special{pa 3197 809}%
\special{pa 3151 855}%
\special{pa 3127 876}%
\special{pa 3103 896}%
\special{pa 3078 916}%
\special{pa 3052 936}%
\special{pa 3026 955}%
\special{pa 3000 973}%
\special{pa 2946 1007}%
\special{pa 2918 1023}%
\special{pa 2890 1038}%
\special{pa 2861 1053}%
\special{pa 2833 1067}%
\special{pa 2804 1080}%
\special{pa 2775 1094}%
\special{pa 2715 1118}%
\special{pa 2655 1140}%
\special{pa 2625 1149}%
\special{pa 2563 1167}%
\special{pa 2501 1181}%
\special{pa 2469 1186}%
\special{pa 2437 1190}%
\special{pa 2405 1193}%
\special{pa 2373 1195}%
\special{pa 2341 1196}%
\special{pa 2309 1196}%
\special{pa 2278 1194}%
\special{pa 2246 1190}%
\special{pa 2214 1183}%
\special{pa 2184 1174}%
\special{pa 2155 1161}%
\special{pa 2127 1144}%
\special{pa 2102 1124}%
\special{pa 2082 1099}%
\special{pa 2069 1070}%
\special{pa 2061 1038}%
\special{pa 2060 1006}%
\special{pa 2065 975}%
\special{pa 2073 944}%
\special{pa 2084 914}%
\special{pa 2099 885}%
\special{pa 2115 857}%
\special{pa 2133 831}%
\special{pa 2152 806}%
\special{pa 2173 781}%
\special{pa 2217 735}%
\special{pa 2240 713}%
\special{pa 2264 691}%
\special{pa 2339 631}%
\special{pa 2365 613}%
\special{pa 2392 595}%
\special{pa 2419 578}%
\special{pa 2447 562}%
\special{pa 2474 546}%
\special{pa 2503 531}%
\special{pa 2559 501}%
\special{pa 2588 487}%
\special{pa 2617 474}%
\special{pa 2647 461}%
\special{pa 2707 439}%
\special{pa 2737 429}%
\special{pa 2799 411}%
\special{pa 2829 403}%
\special{pa 2861 395}%
\special{pa 2892 388}%
\special{pa 2923 383}%
\special{pa 2955 379}%
\special{pa 2987 376}%
\special{pa 3019 374}%
\special{pa 3051 373}%
\special{pa 3083 374}%
\special{pa 3115 376}%
\special{pa 3147 380}%
\special{pa 3178 387}%
\special{pa 3208 397}%
\special{pa 3237 410}%
\special{pa 3265 427}%
\special{pa 3290 447}%
\special{pa 3308 473}%
\special{pa 3319 497}%
\special{fp}%
%
\special{pn 8}%
\special{pa 3690 1730}%
\special{pa 4270 1080}%
\special{fp}%
%
\special{pn 4}%
\special{pa 3100 380}%
\special{pa 2290 1190}%
\special{fp}%
\special{pa 2980 380}%
\special{pa 2190 1170}%
\special{fp}%
\special{pa 2840 400}%
\special{pa 2110 1130}%
\special{fp}%
\special{pa 2660 460}%
\special{pa 2070 1050}%
\special{fp}%
\special{pa 2410 590}%
\special{pa 2090 910}%
\special{fp}%
\special{pa 3200 400}%
\special{pa 2410 1190}%
\special{fp}%
\special{pa 3280 440}%
\special{pa 2550 1170}%
\special{fp}%
\special{pa 3320 520}%
\special{pa 2730 1110}%
\special{fp}%
\special{pa 3300 660}%
\special{pa 2980 980}%
\special{fp}%
%
\special{pn 4}%
\special{pa 4770 630}%
\special{pa 4300 1100}%
\special{fp}%
\special{pa 4670 550}%
\special{pa 4210 1010}%
\special{fp}%
\special{pa 4550 490}%
\special{pa 4120 920}%
\special{fp}%
\special{pa 4430 430}%
\special{pa 4040 820}%
\special{fp}%
\special{pa 4290 390}%
\special{pa 3980 700}%
\special{fp}%
\special{pa 4110 390}%
\special{pa 3950 550}%
\special{fp}%
\special{pa 4870 710}%
\special{pa 4410 1170}%
\special{fp}%
\special{pa 4950 810}%
\special{pa 4530 1230}%
\special{fp}%
\special{pa 5030 910}%
\special{pa 4660 1280}%
\special{fp}%
\special{pa 5090 1030}%
\special{pa 4800 1320}%
\special{fp}%
\put(34.9000,-19.4000){\makebox(0,0)[lb]{$B$}}%
\put(24.8000,-8.7000){\makebox(0,0)[lb]{$D^{n_+}$}}%
\put(43.4000,-9.0000){\makebox(0,0)[lb]{$D^d$}}%
\put(28.9000,-14.6000){\makebox(0,0)[lb]{$l_0$}}%
\put(40.0000,-16.5000){\makebox(0,0)[lb]{$l_0$}}%
%
\special{pn 8}%
\special{pa 2900 1040}%
\special{pa 3390 1750}%
\special{fp}%
%
\special{pn 4}%
\special{sh 1}%
\special{ar 3400 1740 8 8 0 6.2831853}%
%
\special{pn 4}%
\special{sh 1}%
\special{ar 2910 1040 8 8 0 6.2831853}%
\special{sh 1}%
\special{ar 2910 1040 8 8 0 6.2831853}%
\special{sh 1}%
\special{ar 2910 1040 8 8 0 6.2831853}%
%
\special{pn 4}%
\special{sh 1}%
\special{ar 4270 1090 8 8 0 6.2831853}%
\special{sh 1}%
\special{ar 4280 1090 8 8 0 6.2831853}%
\put(29.8000,-11.2000){\makebox(0,0)[lb]{$x_0^{n_+}$}}%
\put(42.1000,-12.1000){\makebox(0,0)[lt]{$x_0^d$}}%
%
\special{pn 4}%
\special{sh 1}%
\special{ar 3690 1730 8 8 0 6.2831853}%
\special{sh 1}%
\special{ar 3690 1730 8 8 0 6.2831853}%
\put(30.9000,-18.5000){\makebox(0,0)[lb]{$y_0^{n_+}$}}%
\put(35.4000,-16.5000){\makebox(0,0)[lb]{$y_0^d$}}%
\end{picture}}%
\caption{Example \ref{lambda1not=0}}
\label{ex2}
 \end{figure}
\end{rei}
\begin{Remark}\label{modify the boundary regularity}
In \cite{I-2-nd terms (2014) C^5 and C^2} for the case $\lambda_1=0$ and in \cite{I-dissipative C^3 boundary and C^2 coefficients(2017)} for the case $D=D^{n_+}$ (or $D^{n_-}$) and $\lambda_0=0$, Ikehata extract the exact form of the top term of the indicator functions for $\tau \to \infty$ under the condition that $\partial D$ is $C^3$ class, $B$ is a ball and $\lambda_1 \in C^2(\partial D)$.  In the next paper, currently in preparation, we will improve the regularities of the boundary $\partial D$ and $\lambda_1$ for Theorem \ref{a result for the case l_0^+ = l_0^-}. 
\end{Remark}
\par
To show Theorem \ref{a result for the case l_0^+ = l_0^-}, we reduce the indicator $I_{\tau}$ to another form $J_{\tau}$, 
which is equivalent to $I_{\tau}$. 
From (\ref{time dependent problem with cavities}) 
and (\ref{Image of partial LT for u}), we can see that $w(x; \tau)$ satisfies
\begin{equation}
\left\{
\begin{array}{ll}(\gamma_0\Delta - \tau^2){\mathcal L}_Tu(x; \tau) + f(x) 
= e^{-{\tau}T}(\partial_tu(T, x)+{\tau}u(T, x)) & \qquad\text{in } \Omega, \\
{\mathcal B^{n}_{\tau}}{\mathcal L}_Tu(x; \tau) 
=\lambda_1(x)e^{-\tau T}u(T,x)  & \qquad\text{on } \partial{D^n}, \\
{\mathcal L}_Tu(x; \tau) =0 & \qquad\text{on } \partial{D^d},
\end{array}
\right.
\label{eq of w in this case}
\end{equation}
where ${\mathcal B_{\tau}^n}w 
= \gamma_0\partial_{\nu_x}w - \lambda(x; \tau)w$ and  
$\lambda(x; \tau) = \lambda_1(x)\tau + \lambda_0(x)$.  
\par
Let $w(x; \tau)$ be the solution of 
\begin{equation}
\left\{
\begin{array}{ll}(\gamma_0\Delta - \tau^2)w(x; \tau) + f(x) = 0 
& \qquad\text{in } \Omega, \\
{\mathcal B^{n}_{\tau}}w(x; \tau) = 0 
& \qquad\text{on } \partial{D^n}, \\
w(x; \tau) = 0 
& \qquad\text{on } \partial{D^d}.
\end{array}
\right.
\label{eq of main part of w}
\end{equation}
Since ${\mathcal L}_Tu(x; \tau) $ and $w$ are $L^2$-weak solution of (\ref{eq of w in this case}) and 
(\ref{eq of main part of w}), respectively, usual elliptic estimate
implies that for some fixed constant $C > 0$ 
$$
\V{{\mathcal L}_Tu(\cdot ; \tau)  - w(\cdot ; \tau)}_{L^2(\Omega)}
\leq Ce^{-{\tau}T}\tau^{-1}
\qquad(\tau >> 1) 
$$
(cf. Appendix C in \cite{Ikehata-Kawashita3}).
The above estimate yields
\begin{equation*}
I_{\tau} = J_{\tau} + O(\tau^{-1}e^{-{\tau}T})  \qquad (\tau \to \infty),
\end{equation*}
where
\begin{equation*}
J_{\tau} = \int_{\Omega}f(x)(w(x; \tau) - v(x; \tau))dx. 
\end{equation*}
\par
Since $w \in H^1_{0, \partial D^d}(\Omega) \subset H^1(\Omega)$ 
and $v \in H^2(\R^3)$ satisfies 
(\ref{reduced background equation}), from integration by parts we obtain 
\begin{align}
\int_{\Omega}fw\,dx &= -\int_{\Omega}({\gamma_0}\Delta-\tau^2)vw\,dx
\nonumber\\&
= \int_{\partial D}\gamma_0\partial_{\nu_x}v\,w\,dS_x 
+\int_{\Omega}(\gamma_0\nabla_x v\cdot\nabla_xw+\tau^2vw)dx. 
\label{integration by parts for 1st term}
\end{align}
A usual elliptic theory implies that the weak solution 
$w \in H^1_{0, \partial D^d}(\Omega)$ of 
(\ref{eq of main part of w}) has the conormal derivative 
$\partial_{\nu_x}w \in H^{-1/2}(\partial\Omega)$ 
in the sense of dual form (see e.g. Lemma 4.3 and Theorem 4.4 in \cite{McLean}). 
From (\ref{eq of main part of w}) it follows that 
\begin{align*}
\int_{\Omega}fv\,dx &= -\int_{\Omega}(\gamma_0\Delta-\tau^2)wv\,dx 
\nonumber
\\&
= \int_{\partial D}
\gamma_0 \partial_{\nu_{x}} w\,v\,dS_x 
+\int_{\Omega}(\gamma_0\nabla_xw\cdot\nabla_xv+\tau^2wv)dx.
\end{align*}
From this equality and (\ref{integration by parts for 1st term}), it follows that
\begin{align*}
J_{\tau} 
= \gamma_0\int_{\partial D}\big(\partial_{\nu_x}v\,w 
- \partial_{\nu_x}w\,v\big)dS_x.
\end{align*}
Since 
$\partial D = \partial{D}^{n}\cup\partial{D}^{d}$, 
$J_{\tau}$ is decomposed into three parts:
\begin{align*}
&J_{\tau} = J_{\tau}^{n_+}+J_{\tau}^{n_-}+J_{\tau}^{d}
=\sum_{\alpha\in \{n_+,n_-,d\}}J_{\tau}^\alpha,  
\end{align*}
where each $J_{\tau}^\alpha$ can be represented by the integral on corresponding boundary as 
\begin{align}
J_{\tau}^{n_\pm}= \int_{\partial{D}^{n_\pm}}{\mathcal B_{\tau}^{n}}v(x; \tau) \,
w(x; \tau)dS_x, 
\quad
J_{\tau}^{d} = - \gamma_0\int_{\partial{D^{d}}}
\partial_{\nu_x}w(x; \tau)\,v(x; \tau)\,dS_x.  
\label{integral on corresponding boundary}
\end{align}
As stated in \cite{IEO2, IEO3 bistatic}, if the cavities are all positive (resp. negative) cavities, the indicator function for large $\tau$ is positive (resp. negative). From these results, we expect that  
$J_{\tau}^{n_+}>0$, $J_{\tau}^{n_-}<0$ and $J_{\tau}^d<0$ for large $\tau$.  
Therefore, even if the parameter $\tau$ is large enough, the sign of the indicator function is not determined. This situation causes difficulties, so it is necessary to pay close attention to the analysis of the indicator function.
In order to analyze the terms (\ref{integral on corresponding boundary}) in detail and obtain the exact form of the top terms of the indicator function, in the next section, we construct an asymptotic solution that approximates the reflection part of the solution $w$ of  (\ref{eq of main part of w}). 

\setcounter{equation}{0}
\section{Construction of an approximate solution}
\label{Approximate solutions of w_0}
For a while, we assume that $D$ has $C^{\ell+2}$ boundary,  
$\lambda_0 \in C^{\ell-2}(\partial D^n)$ and $\lambda_1 \in C^\ell(\partial D^n)$. 
The $\ell \ge 2$ appearing here will be determined later.  
We construct an approximate solution for $w(x; \tau)$ of the form:
$w(x; \tau) = v(x; \tau)\big\vert_{\Omega} + \tilde{w}(x; \tau)$, 
where $\tilde{w}(x; \tau)$ has a kernel representation
$$
\tilde{w}(x; \tau) = \int_{\Omega}\tilde\Psi_\tau(x, y)f(y)\,dy.
$$
Since $v$ has the kernel representation (\ref{free kernel_rep}), 
from (\ref{eq of main part of w}) it follows that the kernel 
$\tilde\Psi_\tau(x, y)$ satisfies
\begin{align}
\left\{
\begin{array}{ll}
(\gamma_0\Delta - \tau^2)\tilde\Psi_\tau(x, y) = 0 
&\qquad\text{in } \Omega, 
\\
{\mathcal B_{\tau}^{n}}
\tilde\Psi_\tau(x, y) 
= 
-{\mathcal B_{\tau}^{n}}\Phi_\tau(x, y) & \qquad\text{on } \partial{D^n}, 
\\ 
\tilde\Psi_\tau(x, y) = -\Phi_\tau(x, y) 
& \qquad\text{on } \partial{D^d}. 
\end{array}
\right.
\label{eq of the kernel tildePsi_tau(x, y) for only cavities}
\end{align}
Here, ${\mathcal B_{\tau}^{n}}\Phi_\tau(x, y)$ can be written by 
\begin{align}
&{\mathcal B_{\tau}^{n}}\Phi_\tau(x, y)
=(-\tau a_0(x,y)+a_1(x,y))e^{-\tau\v{x-y}/\sqrt{\gamma_0}}, 
\label{free expand}
\end{align}
where
\begin{align*} 
a_0 &= \frac{1}{4\pi\v{x - y}}
\Big(\frac{\nu_x\cdot(x- y)}{\sqrt{\gamma_0}\v{x - y}} 
+ \frac{\lambda_1(x)}{{\gamma_0}} \Big), 
\quad
a_1 = \frac{-1}{4\pi\v{x - y}}
\Big(\frac{\nu_x\cdot(x - y)}{\v{x-y}^2}+\frac{\lambda_0(x)}{\gamma_0}\Big). 
\end{align*}
\par
Let us construct an asymptotic solution 
$\Psi_{\tau, N}(x, y)$ of $\tilde\Psi_\tau(x, y)$ satisfying
(\ref{eq of the kernel tildePsi_tau(x, y) for only cavities}) in ${\mathcal U}_4 (\supset
\Gamma)$ of the form:
\begin{equation}
\Psi_{\tau, N}(x, y) 
= e^{-{\tau}\phi(x, y)}\sum_{j = 0}^N(-\tau)^{-j}b_j(x, y)\quad
(x \in {\mathcal U}_4\cap\overline{\Omega}, y \in \overline{B}).
\label{asymptotic_solution_N}
\end{equation}
We also put $b_{-1}(x, y) = 0$.
Acting on (\ref{asymptotic_solution_N}) with the differential operators appearing in
(\ref{eq of the kernel tildePsi_tau(x, y) for only cavities}), we obtain
\begin{align}
({\gamma_0}\Delta& - \tau^2)\Psi_{\tau, N}(x, y) 
= e^{-{\tau}\phi(x, y)}
\Big\{\sum_{j = 0}^N\Big(
(-\tau)^{2-j}({\gamma_0}\v{\nabla_x\phi}^2-1)b_j(x, y)
\nonumber
\\&\hskip-7mm
+ (-\tau)^{1-j}\big(T_{\phi}b_j(x, y)
+{\gamma_0}\Delta b_{j-1}(x, y)\big)\Big)
+(-\tau)^{-N}{\gamma_0}\Delta b_N(x, y)\Big\}
\quad\text{in ${\mathcal U}_4\cap\Omega$, }
\label{equation for the asymptotic solution}
\\
{\mathcal B_{\tau}^{n}}
&\Psi_{\tau, N}(x, y) 
= e^{-{\tau}\phi(x, y)}\Big\{
\sum_{j = 0}^N(-\tau)^{1-j}\big\{
\big((\gamma_0\partial_{\nu_x}{\phi}) + \lambda_1(x)\big)b_j(x, y)
+{\mathcal B_{0}^{n}}b_{j-1}(x,y) \big\} 
\nonumber \\
& \hskip7pc 
+ (-\tau)^{-N}{\mathcal B_{0}^{n}}b_{N}(x, y) \Big\} 
\quad \text{on ${\mathcal U}_4\cap\partial{D^{\alpha}}$ for $\alpha \in \{n_+, n_-\}$,}
\label{equation of boundary conditions for the asymptotic solution}
\end{align}
where 
$T_{\phi} 
= 2(\gamma_0\nabla_x\phi)\cdot\nabla_x + {\gamma_0}\Delta\phi$ and ${\mathcal B_{0}^{n}}w=\gamma_0 \partial_{\nu_x} w-\lambda_0(x) w$. We construct $\phi$ satisfying
\begin{align}
\left\{\begin{array}{ll}
{\gamma_0}\v{\nabla_x\phi}^2 = 1 \quad & \quad\text{in } \Omega\cap{\mathcal U}_4, \\
\phi(x, y) = \v{x - y}/\sqrt{\gamma_0} \quad & \quad\text{on } \partial{D}\cap{\mathcal U}_4, \\
\partial_{\nu_x}\phi(x, y) > 0  \quad & \quad\text{on } \partial{D}\cap{\mathcal U}_4. 
\end{array}
\right.
\label{eikonal eq}
\end{align}
For the amplitude functions $b_j$, we also construct  
\begin{align}
\left\{
\begin{array}{ll}
T_{\phi}b_j(x, y) + {\gamma_0}\Delta b_{j-1}(x, y) = 0 
\quad & \text{in } \Omega\cap{\mathcal U}_4, \\
\big( (\gamma_0 \partial_{\nu_x}{\phi}) + \lambda_1(x) \big) b_j(x,y) 
+{\mathcal B_0^n}b_{j-1}(x,y) = -a_j(x,y)
\quad & \text{on } \partial{D}^n\cap{\mathcal U}_4, 
\\
b_j(x,y) = -\frac{\delta_{0, j}}{4\pi\gamma_0\v{x - y}} 
\quad & \text{on } \partial{D^d}\cap{\mathcal U}_4,
\end{array}
\right.
\label{transport eq}
\end{align}
where 
$\delta_{0, j} = 1$ for $j = 0$, and $\delta_{0, j} = 0$ for $j \geq 1$ and 
$a_j$ ($j=0$, $1$) is given by (\ref{free expand}) and 
$a_j = 0$ for $j \geq 2$.   
Since (\ref{eikonal eq}) and (\ref{transport eq}) are solvable near the boundary, their solutions $\phi$ and $b_j$ exist on ${\mathcal U}_4$, if necessary, by choosing smaller $r_0 > 0$.  
\par
From (\ref{equation for the asymptotic solution}) - (\ref{transport eq}), 
it follows that
\begin{align}
\left\{
\begin{array}{ll}
({\gamma_0}\Delta - \tau^2)\Psi_{\tau, N}(x, y) 
= (-\tau)^{-N}e^{-{\tau}\phi(x, y)}
{\gamma_0}\Delta b_N(x, y)
& \text{in } \Omega\cap{\mathcal U}_4, 
\\
{\mathcal B_{\tau}^{n}}
\Psi_{\tau, N}(x, y) 
= -{\mathcal B_{\tau}^{n}}\Phi_\tau(x, y)  &\\
\hskip3pc+ (-\tau)^{-N}
(\delta_{0,N}a_1(x,y)+{\mathcal B_{0}^{n}}b_{N}(x, y) )e^{-{\tau}\phi(x, y)}
& \text{on } \partial{D}^{n}\cap{\mathcal U}_4, 
\\
\Psi_{\tau, N}(x, y) = -\Phi_\tau(x, y) 
& \text{on } \partial{D}^{d}\cap{\mathcal U}_4,
\end{array}
\right.
\label{equation satisfying the asymptotic equationsD+}
\end{align}
and especially we know that
\begin{equation*}
b_0(x,y)=
\frac{-1}{\gamma_0 \partial_{\nu_x} \phi(x,y)+\lambda_1(x)}a_0 (x,y) \quad 
\text{on $\partial D^\alpha \cap{\mathcal U}^\alpha_3$ for $\alpha \in \{n_+, n_-\}$}.
\end{equation*}
Now, we check the regularity of $\phi$ and $b_j$ ($j=0, 1, \ldots ,N$). Take an open set $\tilde{B}$ satisfying $\overline{B} \subset \tilde{B}$ and $\overline{D}\cap \tilde{B}=\emptyset$. 
We introduce the notation: $\phi \in \tilde{C}^k((\overline{\Omega}\cap{\mathcal U}_4)\times\tilde{B})$, when $\phi(x,y)$ satisfies $\partial_x^{\alpha}\partial_y^{\beta}\phi \in C((\overline{\Omega}\cap{\mathcal U}_4)\times\tilde{B})$ 
in $(x,y) \in (\overline{\Omega}\cap{\mathcal U}_4)\times\tilde{B}$ for $\v{\alpha} \le k$ and $\v{\beta} \ge 0$.  
We construct the solution $\phi(x,y)$ of (\ref{eikonal eq}) due to the way for Eikonal equation (see for example Chap.2 and 3 of Ikawa \cite{{Ikawa book(asymptotic solutions) }}). 
Since $D$ has $C^{\ell+2}$ boundary and $\phi(x, y) \vert_{(\partial{D}\cap{\mathcal U}_4)\times \tilde{B}}= \v{x - y}/\sqrt{\gamma_0}$, it follows that $\phi \in \tilde{C}^{\ell+2}((\overline{\Omega}\cap{\mathcal U}_4)\times\tilde{B})$.   
It means that $\triangle\phi \in \tilde{C}^{\ell}((\overline{\Omega}\cap{\mathcal U}_4)\times\tilde{B})$ and $\gamma_0 (\partial_{\nu_x} \phi)(x,y)+\lambda_1(x) \in \tilde{C}^{\ell}((\partial D^n\cap{\mathcal U}_4)\times\tilde{B})$. 
Therefore, for the solution $b_0$ of the transport equation (\ref{transport eq}) we have $b_0 \in \tilde{C}^\ell((\overline{\Omega}\cap{\mathcal U}_4)\times\tilde{B})$ noting  $b_{-1} = 0$ in (\ref{transport eq}). 
Similarly for $b_1$ in (\ref{transport eq}), we know that $b_1 \in \tilde{C}^{\ell-2}((\overline{\Omega}\cap{\mathcal U}_4)\times\tilde{B})$, since $\Delta b_0 \in \tilde{C}^{\ell-2}((\overline{\Omega}\cap{\mathcal U}_4)\times\tilde{B})$.
If we repeat this argument, we have $b_j \in \tilde{C}^{\ell-2j}((\overline{\Omega}\cap{\mathcal U}_4)\times\tilde{B})$ for $j=1, \ldots ,N$.  Set $\ell=2N+2$, then $\Delta b_j \in \tilde{C}^{2(N-j)}((\overline{\Omega}\cap{\mathcal U}_4)\times\tilde{B}) \subset C^0((\overline{\Omega}\cap{\mathcal U}_4)\times\overline{B})$ ($j=0, 1, \ldots ,N$), that is consistent with 
(\ref{equation satisfying the asymptotic equationsD+}).  
\par
From the above, we use the following terminology for convenience.   
\begin{Def}\label{N-class}
We call $(\partial D, \lambda)$ is in $m$-class if $\partial D$ is $C^{2m+4}$, $\lambda_1 \in C^{2m+2}(\partial D^n)$ and $\lambda_0 \in C^{2m}(\partial D^n)$ for integer $m\ge 1$, and $(\partial D, \lambda)$ is in $0$-class(i.e., $m=0$) if $\partial D$ is  $C^{4}$, $\lambda_1 \in C^{2}(\partial D^n)$ and $\lambda_0 \in L^{\infty}(\partial D^n)$. 
\end{Def}
Hereafter, we proceed with the discussion assuming that $(\partial D, \lambda)$ is in $N$-class, where $\phi \in \tilde{C}^{2N+4}((\overline{\Omega}\cap{\mathcal U}_4)\times\tilde{B})$ and $b_j \in \tilde{C}^{2(N-j)+2}((\overline{\Omega}\cap{\mathcal U}_4)\times\tilde{B})$, $j=1, \ldots, N$. 
When $N=0$, we have $b_0 \in \tilde{C}^{2}((\overline{\Omega}\cap{\mathcal U}_4)\times\tilde{B})$ even for $\lambda_0 \in L^{\infty}(\partial D^n)$ since $b_{-1}=0$ and the definition of $a_0$. 
The minimum order $N$ required to obtain the result will be determined later. 
\par
Choose a cutoff function $\chi \in C^\infty_0(\R^3)$ with $\chi(x) = 1$ for $x \in {\mathcal U}_1$, ${\rm supp}\chi \subset {\mathcal U}_2$. We put 
\begin{equation}
\Psi_{\tau, N}^r(x, y) = \tilde\Psi_\tau(x, y) 
- \chi(x)\Psi_{\tau, N}(x, y). 
\label{definition of the remainder term}
\end{equation}
From (\ref{eq of the kernel tildePsi_tau(x, y) for only cavities}) and (\ref{equation satisfying the asymptotic equationsD+}), we can see that $\Psi_{\tau, N}^r(x, y)$ satisfies 
\begin{equation}
\left\{
\begin{array}{ll}
({\gamma_0}\Delta - \tau^2)\Psi_{\tau, N}^r(x, y) = F_{\tau, N}(x, y)
&\qquad\text{in } \Omega, \\
{\mathcal B_{\tau}^n}\Psi_{\tau, N}^r(x, y) 
= G_{\tau, N}^n(x, y)
& \qquad\text{on } \partial{D}^n,  
\\ 
\Psi_{\tau, N}^r(x,y) = H^d_{\tau}(x,y)
& \qquad\text{on } \partial{D}^{d},
\end{array}
\right.
\label{equations of the remainder term (for cavities only)}
\end{equation}
where $F_{\tau, N}$ and $G_{\tau, N}^n$ are given by
\begin{align}
F_{\tau, N}(x, y) 
&= -\Big\{(-\tau)^{-N}\chi(x)e^{-{\tau}\phi(x, y)}
{\gamma_0}\Delta b_N(x, y) 
+ \gamma_0[\Delta, \chi]\Psi_{\tau, N}(x, y)\Big\}, 
\label{definition of F_N}
\\
G_{\tau, N}^n(x, y) &= (\chi(x)-1)
{\mathcal B_{\tau}^n}\Phi_{\tau}(x, y)
\nonumber
\\& \hskip-4pc
-\Big\{(-\tau)^{-N}e^{-{\tau}\phi(x, y)}\chi(x)
(\delta_{0,N}a_1(x,y)+{\mathcal B_{0}^n}b_{N}(x, y))
+\gamma_0\partial_{\nu_x}\chi(x)\Psi_{\tau, N}(x, y)\Big\},
\label{definition of G_N}
\\
H^d_{\tau}(x,y) &=(\chi(x)-1)\Phi_{\tau}(x,y).
\label{definition of H_tau}
\end{align}
\par
Using $\Phi_\tau(x, y)$ and $\Psi_{\tau, N}(x, y)$, we can express $w(x; \tau)$ as  
\begin{align*}
w(x;\tau)&=\int_B \big(\Phi_{\tau}(x,{y})
+\tilde{\Psi}_{\tau}(x,{y})\big) f({y})\, d{y} 
= w_{N}(x; \tau) + w_{N}^r(x; \tau), 
\end{align*}
where
\begin{align*}
w_{N}(x; \tau) &= \int_B \big(\Phi_{\tau}(x,{y})
+ \chi(x)\Psi_{\tau, N}(x, {y})\big)f({y})\, d{y}, 
\\
\quad
w_{N}^r(x; \tau) &= \int_B {\Psi}_{\tau,N}^r(x,{y}) 
f({y})\, d{y}. 
\end{align*}
For $\alpha\in\{n_+, n_-\}$, we divide $J_{\tau}^{\alpha}$ into two parts: 
$J_{\tau}^\alpha= J_{\tau, N}^\alpha + J_{\tau, N}^{\alpha,r}$, where
\begin{align*}
J_{\tau, N}^\alpha &= \int_{\partial{D^\alpha}}{\mathcal B_{\tau}^n}v(x; \tau) 
w_{N}(x; \tau)dS_x, 
\quad
J_{\tau, N}^{\alpha,r} = \int_{\partial{D^\alpha}}{\mathcal B_{\tau}^n}v(x; \tau) 
w_{N}^r(x; \tau)dS_x. 
\end{align*}
Similarly, we divide 
$J_{\tau}^{d}$ into $J_{\tau}^{d} = J_{\tau, N}^{d} + J_{\tau, N}^{d,r}$, where
\begin{align*}
J_{\tau, N}^{d} &= - \gamma_0\int_{\partial{D^{d}}}
\partial_{\nu_x}w_{N}(x; \tau)v(x; \tau)dS_x, 
\quad
J_{\tau, N}^{d,r} &= - \gamma_0\int_{\partial{D^{d}}}
\partial_{\nu_x}w_{N}^r(x; \tau)v(x; \tau)dS_x. 
\end{align*}
In section \ref{Estimates of the remainder terms}, we will prove the following lemma:
\begin{Lemma}\label{estimates of the remainder terms for the indicator function}
Assume that $(\partial D, \lambda)$ is in $N$-class. Then, there exists a constant $C > 0$ such that
\begin{align*}
\sum_{\alpha \in \{n_+, n_-, d\}}\v{J_{\tau, N}^{\alpha,r}}
&\leq C{\tau}^{-N-5}e^{-2{\tau}l_0/\sqrt{\gamma_0}}
\qquad (\tau \geq 1).
\end{align*}
\end{Lemma}
For any $\alpha \in \{n_+, n_-, d\}$, the term $J_{\tau,N}^{\alpha}$ corresponds  
to the main part of $J_{\tau}^{\alpha}$. 
Thus, the main term of $J_{\tau,N}^{\alpha}$ determines the structure of the main term of $J_{\tau}^{\alpha}$. In Section \ref{Top term structure}, the asymptotic behavior of $J_{\tau,N}^{\alpha}$ and the structure of the top terms are investigated using the Laplace method below:
\begin{Lemma}\label{a usual Laplace method} 
Let $U$ be an arbitrary open set of ${\R}^n$.  
Assume that $h \in C^{2,\beta_0}(\overline U)$ with some $\beta_0 > 0$ and  
$h(x) > h(x_0)$ for all $x\in\overline U\setminus\{x_0\}$ 
at a point $x_0\in U$, and $\text{det}\,(\text{Hess}\,(h)(x_0))>0$. 
Then for given $\varphi\in C^{0,\beta_0}(\overline U)$ 
it holds that
$$\displaystyle
\int_U e^{-\tau h(x)}\varphi(x)dx
=\frac{e^{-\tau h(x_0)}}{\sqrt{\text{det}\,(\text{Hess}\,(h)(x_0))}}
\left(\frac{2\pi}{\tau}\right)^{n/2}
\left(\varphi(x_0)+\Vert\varphi\Vert_{C^{0,\beta_0}(\overline U)}
O(\tau^{-\beta_0/2})\right).
$$
Moreover there exists a positive constant $C$ such that, 
for all $\tau \geq 1$
$$\displaystyle
\left\vert\int_Ue^{-\tau h(x)}\varphi(x)dx\right\vert\le
\frac{C e^{-\tau h(x_0)}}{\tau^{n/2}}
\Vert\varphi\Vert_{C(\overline U)}.
$$
\end{Lemma}
For a proof, see e.g. Appendix A in \cite{Ikehata-Kawashita3}. 
\par

\setcounter{equation}{0}
\section{Asymptotic form of the main terms of $I_{\tau}$} 
\label{Top term structure}

\par
Once, we assume that $(\partial D, \lambda)$ is in $N$-class for the integer $N \ge 0$ and we admit Lemma \ref{estimates of the remainder terms for the indicator function}. Let us consider $J_{\tau, N}^\alpha$ for $\alpha \in\{n_+, n_-\}$:
\begin{align*}
&J_{\tau, N}^\alpha
=\int_{\partial{D^\alpha}}
{\mathcal B_{\tau}^n}v(x; \tau) 
w_{N}(x; \tau)dS_x
\nonumber
\\
&=\int_{B\times B}
dyd\tilde{y}
f(y)f(\tilde{y})\int_{\partial D^\alpha}
\big(\Phi_{\tau}(x, \tilde{y})+\chi(x){\Psi}_{\tau, N}(x,\tilde{y})\big) 
{\mathcal B_{\tau}^n} \Phi_{\tau}(x,y)\,dS_x.
\nonumber
\end{align*}
Since $\phi(x,y)=\vert x-y \vert/\sqrt{\gamma_0}$ on $\partial D \cap {\mathcal U}_3$ and (\ref{free expand}), we have the expansion: $J_{\tau,N}^\alpha=
\sum_{k=-1}^N \tau^{-k}K_{\tau,-k}^\alpha$, where the top term $K_{\tau,1}^\alpha$ is described as 
\begin{align*}
K_{\tau,1}^\alpha&=\int_{\partial D^\alpha \times B \times B }
f(y)f(\tilde{y})e^{-\tau L(x,y, \tilde{y})/\sqrt{\gamma_0}} 
\kappa_{1}^\alpha(x,y,\tilde{y})
dS_x dy d\tilde{y}, 
\\
\nonumber
&\kappa_{1}^\alpha(x,y,\tilde{y}) =\eta^{\alpha} (x, y)\tilde{\eta}^{\alpha}(x, \tilde{y}) 
\quad((x, y, \tilde{y}) \in (\partial{D}^{\alpha}{\cap}{\mathcal U}_3) \times \overline{B} \times \overline{B}), 
\\
&\eta^{\alpha}(x,y) =-a_0(x,y),\quad
\tilde{\eta}^{\alpha}(x,\tilde{y})=
\frac{1}{4\pi \gamma_0\v{x-\tilde{y}}}+\chi(x)b_0(x,\tilde{y}). 
\end{align*}
The lower order terms $K_{\tau,-k}^\alpha$ ($k=0, 1, 2, \ldots, N$) is described as 
\[
K_{\tau,-k}^\alpha=\int_{\partial D^\alpha \times B \times B }
f(y)f(\tilde{y})e^{-\tau L(x,y, \tilde{y})/\sqrt{\gamma_0}}
\kappa_{-k}^\alpha(x, y,\tilde{y}) dS_x dy d\tilde{y},
\]
where the kernels $\kappa_{-k}^\alpha(x, y,\tilde{y})$ for $(x, y, \tilde{y}) \in (\partial{D}^{\alpha}{\cap}{\mathcal U}_3) \times \overline{B} \times \overline{B}$ are determined depending on the order $N$. For $N=0$, we have 
\[\kappa_{0}^\alpha(x,y,\tilde{y}) 
=a_1(x,y)\left(\frac{1}{4\pi \gamma_0\v{x-\tilde{y}}}
+\chi(x)b_0(x,\tilde{y})\right). \]
For $N=1$, we have 
\begin{align*}
&\kappa_{0}^\alpha(x,y,\tilde{y}) 
=a_1(x,y)\left(\frac{1}{4\pi \gamma_0\v{x-\tilde{y}}}
+\chi(x)b_0(x,\tilde{y})\right)
+\chi(x)a_0(x,y)b_1(x,\tilde{y}),  \enskip 
\\&\kappa_{-1}^\alpha(x,y,\tilde{y})
=(-1)^N\chi(x) a_1(x,y)b_{1}(x,\tilde{y}).   
\end{align*}
For $N \geq 2$, $\kappa_{0}^\alpha(x,y,\tilde{y})$ is the same as $N=1$, and for $k=1, 2, \ldots, N$, 
\begin{align*}
&\kappa_{-k}^\alpha(x,y,\tilde{y})
=(-1)^k\chi(x)\left(a_0(x,y)b_{k+1}(x,\tilde{y})
+a_1(x,y)b_{k}(x,\tilde{y}) \right), \\
&\kappa_{-N}^\alpha(x,y,\tilde{y})
=(-1)^N\chi(x) a_1(x,y)b_{N}(x,\tilde{y}). 
\end{align*}

Next we consider $J_{\tau,N}^{d}$:
\begin{align}
J_{\tau,N}^{d}
&=- \gamma_0\int_{\partial{D^{d}}}\partial_{\nu_x}
w_{N}(x; \tau)v(x; \tau)dS_x
\nonumber
\\
&=-\int_{\partial D^{d} \times B \times B}
dS_x dy d\tilde{y} f(y)f(\tilde{y}) \gamma_0\partial_{\nu_x}
\big\{\Phi_{\tau}(x,\tilde{y})+\chi(x)\Psi_{\tau,N}^{d}(x,\tilde{y}) \big\}
\Phi_{\tau}(x,y).    
\nonumber
\end{align}
Thus, in similar way to the case for $J_{\tau,N}^{n_\pm}$, 
we have 
$J_{\tau,N}^{d}=
-\sum_{k=-1}^N \tau^{-k}K^{d}_{\tau,-k}$, 
where the term $K_{\tau,1}^{d}$ is described as
\begin{align*}
K_{\tau,1}^{d} =&\int_{\partial D^{d} \times B \times B}
f(y)f(\tilde{y})e^{-\tau L(x,y,\tilde{y})/\sqrt{\gamma_0}}
\kappa_{1}^d(x, y,\tilde{y}) 
dS_x dy d\tilde{y},  
\\
&\kappa_{1}^d(x, y,\tilde{y})
= \eta^{d}(x,y) \tilde{\eta}^{d}(x,\tilde{y})
\quad ((x, y, \tilde{y}) \in (\partial{D}^{d}{\cap}{\mathcal U}_3) \times \overline{B} \times \overline{B}), 
\\
&\eta^{d}(x,y)
=\frac{1}{4\pi\gamma_0\v{x-y}}, \\
&\tilde{\eta}^{d}(x,\tilde{y})
=\frac{1}{4\pi\sqrt{\gamma_0}\v{x-\tilde{y}}}
\Big\{
-\nu_x \cdot\frac{x-\tilde{y}}{\v{x-\tilde{y}}}
+\chi(x)
\sqrt{\gamma_0}\partial_{\nu_x}\phi(x,\tilde{y})
\Big\}, 
\end{align*} 
and the terms $K_{\tau,-k}^{d}$ ($k=0, 1, 2, \ldots, N$) are 
described as 
\begin{align*}
&K_{\tau,-k}^{d}=\int_{\partial D^{d} \times B \times B }
f(y)f(\tilde{y})e^{-\tau L(x,y, \tilde{y})/\sqrt{\gamma_0}}
\kappa_{-k}^{d}(x, y,\tilde{y}) dS_x dy d\tilde{y},
\\
&\kappa_{0}^{d}(x,y,\tilde{y}) = 
\frac{-1}{4\pi\v{x-y}}
\Big(
\frac{\nu_x\cdot \frac{x-\tilde{y}}
{\v{x-\tilde{y}}}}{4\pi\gamma_0\v{x-\tilde{y}}^2}+\partial_{\nu_x}\chi(x)\frac{1}{4\pi\gamma_0\v{x-\tilde{y}}}\Big)
+\frac{\chi(x)}{4\pi\v{x-y}}\partial_{\nu_x}b_0(x,\tilde{y}), \\
&\kappa_{-k}^d(x,y,\tilde{y})
=\frac{(-1)^{k}}{4\pi\v{x-y}}\chi(x)\partial_{\nu_x}b_k(x,\tilde{y}) 
\quad \text{for $k=1, 2, \ldots, N$ as $N \ge 1$}\\
&\hskip10pc\quad((x, y, \tilde{y}) \in (\partial{D}^{d}{\cap}{\mathcal U}_3) \times \overline{B} \times \overline{B}).  
\end{align*}
\par
By using these expansions of $J_{\tau, N}^\alpha$ for $\alpha \in\{n_+, n_-,d\}$, we find the exact form of the top term of the indicator function and prove Theorem \ref{a result for the case l_0^+ = l_0^-_general}.  

\par

{\it Proof of Theorem \ref{a result for the case l_0^+ = l_0^-_general}.}  
To evaluate the term $K_{\tau,1}^\alpha$ for $\alpha \in \{n_+, n_-, d\}$, 
it is important to note that 
\begin{equation}
L(x, y, \tilde{y}) \geq 2(l_0+c_0 )
\quad((x, y, \tilde{y}) \in (\partial{D} \times \overline{B} \times \overline{B})\setminus({\mathcal U}_1\times {\mathcal V}_1 \times {\mathcal V}_1) )
\label{estimates on the outside of E^alpha 2}
\end{equation}
holds from (\ref{estimates on the outside of E^alpha}).  
Hence, if we set 
\begin{align*}
K_{\tau,1}^{0,\alpha}(x_0^\alpha, y_0^\alpha) = \int_{\partial D^\alpha \cap B_{r_0}(x_0^\alpha)}& 
\Big( \int_{\{B \cap B_{r_0}(y_0^\alpha)\}^2 }
f(y)f(\tilde{y})e^{-\tau L(x,y,\tilde{y})/\sqrt{\gamma_0}} 
\\&
\eta^{\alpha} (x, y) \tilde{\eta}^{\alpha}(x, \tilde{y})
dyd\tilde{y} \, \Big)dS_x
\quad((x_0^\alpha, y_0^\alpha) \in E^\alpha_0),
\end{align*}
from (\ref{estimates on the outside of E^alpha 2}) it follows that
$$
K_{\tau,1}^\alpha = \sum_{j = 1}^{M^\alpha}K_{\tau,1}^{0,\alpha}(x_j^\alpha, y_j^\alpha) + O(e^{-\frac{2\tau}{\sqrt{\gamma_0}}(l_0+c_0)})\quad(\tau \to \infty).
$$
%
Take any point $(x_0^{\alpha}, y_0^\alpha, \tilde{y}_0^\alpha) \in E^\alpha$, 
Lemma \ref{properties of the set of the minimizers} says that
$\tilde{y}_0^\alpha=y_0^\alpha$ and $x_0^\alpha \in \partial D^{\alpha}$ and 
$y_0^\alpha \in \partial B$. 
Using the local coordinate for $y^b \in \partial{B}\cap{B_{4r_0}(y_0^\alpha)}$ introduced in section \ref{Non-degenerate condition}: 
\begin{align*}
y^b=b^\alpha(u)=x_0^\alpha +u_1e_1+u_2 e_2 +( l_0 + h^\alpha(u))\nu_{x_0^\alpha}
\quad(u \in U_{y_0^\alpha}), 
\end{align*} 
we denote $y \in B{\cap}B_{4r_0}(y_0^\alpha)$ by 
\begin{align}
y=y^\alpha (u, u_3)&=x_0^\alpha +u_1e_1+u_2 e_2 +(l_0 + h^\alpha(u)+u_3)\nu_{x_0^\alpha}\quad(u \in U_{y_0^\alpha}, u_3 \geq 0).   
\label{localB}
\end{align}
Set $\tilde{U}_{y_0^\alpha} = \{ u \in U_{y_0^\alpha} \mid b^\alpha(u) \in \partial{B}{\cap}B_{r_0}(y_0^\alpha)\}$ and $t^\alpha_0(u) = \sup\{ u_3 \mid y^\alpha(u,u_3) \in B{\cap}B_{3r_0}(y_0^\alpha)\}$ $(u \in \tilde{U}_{y_0^\alpha})$.  
By choosing $r_0 > 0$ smaller if necessary, we can assume that $r_0 \leq t^\alpha_0(u) \leq 3r_0$ $(u \in \tilde{U}_{y_0^\alpha})$.  
For another point $\tilde{y}^b \in \partial{B}\cap{B_{r_0}(y_0^\alpha)}$, $\tilde{y}^b$ is denoted by $\tilde{y}^b = b^\alpha(\tilde{u})$ ($\tilde{u} \in \tilde{U}_{y_0^\alpha}$).  
In the same way, for $x \in \partial{D}{\cap}B_{r_0}(x_0^\alpha)$, 
$x$ is denoted by 
$$
x=s^\alpha(\sigma)=x_0^\alpha +\sigma_1e_1+\sigma_2 e_2 - g^\alpha(\sigma_1, \sigma_2)\nu_{x_0^\alpha}
\quad(\sigma = (\sigma_1, \sigma_2) 
\in \tilde{U}_{x_0^\alpha} (\subset U_{x_0^\alpha})),
$$
where
$\tilde{U}_{x_0^\alpha} = \{\sigma \in U_{x_0^\alpha} \mid s^\alpha(\sigma) \in \partial{D}{\cap}B_{r_0}(x_0^\alpha)\}$.  

\par

To obtain the top term of the expansion, it is enough to estimate $K_{\tau,1}^{0,\alpha}(x_0^\alpha, y_0^\alpha)$ for $(x_0^\alpha, y_0^\alpha) \in E^\alpha_0$.
Since we assume that $(\partial D, \lambda)$ is in $N$-class and $N \ge 0$, we have $\phi \in \tilde{C}^{4}((\overline{\Omega}\cap{\mathcal U}_4)\times\tilde{B})$ and $b_0 \in \tilde{C}^{2}((\overline{\Omega}\cap{\mathcal U}_4)\times\tilde{B})$ even for $N=0$. 
Noting that 
\begin{align}
\vert \eta^{\alpha} (x, y) \tilde{\eta}^{\alpha}(x, \tilde{y})
-\eta^{\alpha} (x, y^b) \tilde{\eta}^{\alpha}(x, \tilde{y}^b)\vert 
\le C(\v{y-y^b}+ \v{\tilde{y}-\tilde{y}^b}) 
\le \tilde{C}(u_3+\tilde{u}_3), \label{difference}
\end{align} 
we decompose $K_{\tau,1}^{0,\alpha}(x_0^\alpha, y_0^\alpha)$ into 
$K_{\tau,1}^{0,\alpha}(x_0^\alpha, y_0^\alpha)=K_{\tau,1}^{00,\alpha}(x_0^\alpha, y_0^\alpha)+\tilde{K}_{\tau,1}^{0,\alpha}(x_0^\alpha, y_0^\alpha)$, 
where 
\begin{align*}
K_{\tau,1}^{00,\alpha}(x_0^\alpha, y_0^\alpha)=\int_{\partial D^\alpha \cap B_{r_0}(x_0^\alpha)} 
&\left(
\int_{B \cap B_{r_0}(y_0^\alpha)}
f(y)e^{-\tau \v{x-y}/\sqrt{\gamma_0}} 
\eta^{\alpha} (x, y^b) dy \right)
\\
&\qquad\left(\int_{B \cap B_{r_0}(y_0^\alpha)}
f(\tilde{y})e^{-\tau \v{x-\tilde{y}}/\sqrt{\gamma_0}}
\tilde{\eta}^{\alpha}(x, \tilde{y}^b) d\tilde{y} \right)\,dS_x, 
\end{align*}
\begin{align*}
\tilde{K}_{\tau,1}^{0,\alpha}(x_0^\alpha, y_0^\alpha)=\int_{\partial D \cap B_{r_0}(x_0^\alpha)} 
\Big\{ &\int_{\{B \cap B_{r_0}(y_0^\alpha)\}^{2} }
f(y)f(\tilde{y})e^{-\tau L(x,y,\tilde{y})/\sqrt{\gamma_0}} \\
&\big(\eta^{\alpha} (x, y) \tilde{\eta}^{\alpha}(x, \tilde{y})
-\eta^{\alpha} (x, y^b) \tilde{\eta}^{\alpha}(x, \tilde{y}^b)\big)
dyd\tilde{y} \, \Big\}dS_x. 
\end{align*}
In $K_{\tau,1}^{00,\alpha}(x_0^\alpha, y_0^\alpha)$,  
regarding (\ref{localB}), we have 
\begin{align*}
&\int_{B \cap B_{r_0}(y_0^\alpha)}
f(y)e^{-\tau \v{x-y}/\sqrt{\gamma_0}} 
\eta^{\alpha}(x, y^b) dy \\
&=\int_{\tilde{U}_{y_0^\alpha}} \eta^{\alpha} (x, b^\alpha({u}))
\left( \int_0^{t_0^\alpha(u)} f(y^\alpha(u,u_3)) 
e^{-\tau \v{x-y^\alpha(u,u_3)}/\sqrt{\gamma_0}} \,d u_3 \right)
\,J_{\partial B}^\alpha(u)du,
\end{align*} 
where
$J_{\partial B}^\alpha(u)=\sqrt{1+(h^\alpha_{u_1})^2+(h^\alpha_{u_2})^2}$.  
Set
$X^{\alpha}(u,u_3,\sigma)=\{2(l_0 + k^\alpha(u,\sigma))+u_3\}u_3$.  
For $x=s^\alpha(\sigma)$ on $\partial D\cap B_{r_0}(x_0^\alpha)$, the expansion  
\begin{align*}
\v{x-y^\alpha(u,u_3)}&=\v{s^\alpha(\sigma)-y^\alpha(u,u_3)}=
\sqrt{\v{s^\alpha(\sigma)-b^\alpha(u)}^2+X^{\alpha}(u,u_3,\sigma)} \\
&=\v{s^\alpha(\sigma)-b^\alpha(u)}
\left(1+\frac{X^{\alpha}(u,u_3,\sigma)}{2\v{s^\alpha(\sigma)-b^\alpha(u)}^2}
-\frac{X^{\alpha}(u,u_3,\sigma)^2}{8\v{s^\alpha(\sigma)-b^\alpha(u)}^4}\right.
\\&\hskip70mm
\left. +\frac{X^{\alpha}(u,u_3,\sigma)^3}{16\v{s^\alpha(\sigma)-b^\alpha(u)}^6}-\cdots \right)
\end{align*}
uniformly converges for $\v{X^{\alpha}(u,u_3,\sigma)} \le  \v{s^\alpha(\sigma)-b^\alpha(u)}^2/2$. 
Since $\v{s^\alpha(\sigma)-b^\alpha(u)} \ge l_0$ and 
we can take $r_0>0$ small enough if necessary, the convergence radius is positive. 
We denote the above terms by
\begin{align*}
\v{s^\alpha(\sigma)-y^\alpha(u,u_3)}
=\v{s^\alpha(\sigma)-b^\alpha(u)}+\frac{ l_0 + k^\alpha(u,\sigma)}{\v{s^\alpha(\sigma)-b^\alpha(u)}}u_3 
+B^\alpha(u,u_3,\sigma)u_3^2.
\end{align*}
If we put $A^\alpha(u, \sigma)=\frac{l_0 + k^\alpha(u,\sigma)}{\v{s^\alpha(\sigma)-b^\alpha(u)}}$ and 
$\tilde{f}^\alpha(u,u_3)=f(y^\alpha(u,u_3))$, 
then we have
\begin{align*}
&\int_0^{t_0^\alpha(u)} f(y^\alpha(u,u_3)) 
e^{-\tau \v{x-y^\alpha(u,u_3)}/\sqrt{\gamma_0}} \,d u_3\\
&= e^{-\frac{\tau}{\sqrt{\gamma_0}} \v{s^\alpha(\sigma)-b^\alpha(u)}}
\int_0^{t_0^\alpha(u)}\tilde{f}^\alpha(u,u_3)
e^{-\frac{\tau}{\sqrt{\gamma_0}} A^\alpha(u,\sigma)u_3}
e^{-\frac{\tau}{\sqrt{\gamma_0}} B^\alpha(u,u_3,\sigma)u_3^2}\, du_3
\\
&= e^{-\frac{\tau}{\sqrt{\gamma_0}}\v{s^\alpha(\sigma)-b^\alpha(u)}}
\left\{ \frac{\tilde{f}^\alpha(u,0)}{\frac{\tau}{\sqrt{\gamma_0}} A^\alpha(u,\sigma)}  
+ \frac{e^{-\frac{\tau}{\sqrt{\gamma_0}} A^\alpha(u,\sigma)t_0^\alpha(u)}}
{-\frac{\tau}{\sqrt{\gamma_0}} A^\alpha(u,\sigma)} 
e^{-\frac{\tau}{\sqrt{\gamma_0}} B^\alpha(u,t_0^\alpha(u),\sigma)t_0^\alpha(u)^2} 
\tilde{f}^\alpha(u,t_0^\alpha(u))\right. 
\\
&\qquad +\frac{1}{\frac{\tau}{\sqrt{\gamma_0}} A^\alpha(u,\sigma)}
\int_0^{t_0^\alpha(u)} e^{-\frac{\tau}{\sqrt{\gamma_0}} A^\alpha(u,\sigma)u_3}
\left(\frac{\partial}{\partial u_3}\tilde{f}^\alpha(u,u_3) \right) 
e^{-\frac{\tau}{\sqrt{\gamma_0}} B^\alpha(u,u_3,\sigma)u_3^2 }\,du_3 \\
&\qquad -\frac{1}{A^\alpha(u,\sigma)}
\int_0^{t_0^\alpha(u)}e^{-\frac{\tau}{\sqrt{\gamma_0}} A^\alpha(u,\sigma)u_3}
\tilde{f}^\alpha(u,u_3)
e^{-\frac{\tau}{\sqrt{\gamma_0}} B^{\alpha}(u,u_3,\sigma)u_3^2}
\\
&\left. \qquad \qquad \qquad \qquad \qquad \qquad 
\left( \frac{\partial}{\partial u_3}B^\alpha(u,u_3,\sigma) u_3
+2B^{\alpha}(u,u_3,\sigma)\right) 
\,u_3\,du_3 \right\}\\
&=: e^{-\frac{\tau}{\sqrt{\gamma_0}} \v{s^\alpha(\sigma)-b^\alpha(u)}}
\left\{  
\frac{\tilde{f}^\alpha(u,0)}{\frac{\tau}{\sqrt{\gamma_0}} A^\alpha(u,\sigma)}
+I_1^{1,\alpha} +I_2^{1,\alpha} + I_3^{1,\alpha}
\right\}.
\end{align*}
 If necessary, we choose $r_0 > 0$ for the coverings 
 ${\mathcal U}^\alpha_3$ of $E_0^\alpha$ small enough again.
Then, there exists a $\delta_0 >0$ such that 
$\{A^\alpha(u, \sigma)+B^\alpha(u,t_0^\alpha(u),\sigma)t_0^\alpha(u)\}t_0^\alpha(u) \ge \delta_0$ ($u \in \tilde{U}_{y_0^\alpha}$)  
and we have $I_1^{1,\alpha}=O(e^{-\frac{\tau}{\sqrt{\gamma_0}}\delta_0})$.  
Putting $s=\frac{\tau}{\sqrt{\gamma_0}}A^\alpha u_3$, we have 
\begin{align}
\v{I_2^{1,\alpha}}&\le C \int_0^{\infty} e^{-s}\,ds 
\frac{{\gamma_0}}{\Big(\tau A^\alpha(u,\sigma)\Big)^2}=O(\tau^{-2}),  
\label{I_2}\\
\v{I_3^{1,\alpha}}&\le  C \int_0^{\infty} e^{-s}s\,ds 
\frac{{\gamma_0}}{\Big(\tau A^\alpha(u,\sigma)\Big)^2}=O(\tau^{-2}). \label{I_3}
\end{align}
Thus, it follows that 
\begin{align*}
&\int_{B \cap B_{r_0}(y_0^\alpha)}
f(y)e^{-\tau \v{x-y}/\sqrt{\gamma_0}} 
\eta^{\alpha} (x, y^b) dy \\
=&\frac{\sqrt{\gamma_0}}{\tau}\int_{\tilde{U}_{y_0^\alpha}}
\eta^{\alpha} (s^\alpha(\sigma), b^\alpha(u)) 
e^{-\frac{\tau}{\sqrt{\gamma_0}} \v{s^\alpha(\sigma)-b^\alpha(u)}}
\left(
\frac{{f}(b^\alpha(u))}{ A^\alpha(u,\sigma)}+A_{\tau}^{1,\alpha}(u,\sigma)
\right)J_{\partial B}^\alpha(u)\, du,
\end{align*}
where $A_{\tau}^{1,\alpha}:=\frac{\tau}{\sqrt{\gamma_0}}\{I_1^{1,\alpha} +I_2^{1,\alpha} +I_3^{1,\alpha}\}$ satisfies that 
$A_{\tau}^{1,\alpha}=O(\tau^{-1})$ as $\tau \to \infty$.  
In the same way as above, we have
\begin{align*}
&\int_{B \cap B_{r_0}(y_0^\alpha)}
f(\tilde{y})e^{-\tau \v{x-\tilde{y}}/\sqrt{\gamma_0}}
\tilde{\eta}^{\alpha}(x, \tilde{y}^b) d\tilde{y} \\
=&\frac{\sqrt{\gamma_0}}{\tau}\int_{\tilde{U}_{y_0^\alpha}}
\tilde{\eta}^{\alpha} (s^\alpha(\sigma), b^\alpha(\tilde{u})) 
e^{-\frac{\tau}{\sqrt{\gamma_0}} \v{s^\alpha(\sigma)-b^\alpha(\tilde{u})}}
\left(
\frac{{f}(b^\alpha(\tilde{u}))}{ A^\alpha(\tilde{u},\sigma)}+A_{\tau}^{2,\alpha}(\tilde{u},\sigma)
\right)J_{\partial B}^\alpha(\tilde{u})\, d\tilde{u},
\end{align*}
where $A_{\tau}^{2,\alpha}=O(\tau^{-1})$ as $\tau \to \infty$.
\par
For the top term $K_{\tau,1}^{00,\alpha}(x_0^\alpha, y_0^\alpha)$ of $K_{\tau,1}^{0,\alpha}(x_0^\alpha, y_0^\alpha)$, there exists $B_{\tau}^\alpha(u, \tilde{u}, \sigma)$ 
such that
\begin{align*}
K_{\tau,1}^{00,\alpha}(x_0^\alpha, y_0^\alpha) &=
\frac{\gamma_0}{{\tau}^2} 
\int_{\tilde{U}_{y_0^\alpha}\times \tilde{U}_{y_0^\alpha} \times \tilde{U}_{x_0^\alpha} }
e^{-\frac{\tau}{\sqrt{\gamma_0}}\tilde{L}^\alpha(\sigma, u, \tilde{u})}
\Bigm( f(b^\alpha(u))f(b^\alpha(\tilde{u}))
\\
&
\eta^{\alpha} (s^\alpha(\sigma), b^\alpha(u))\tilde{\eta}^{\alpha}(s^\alpha(\sigma), b^\alpha(\tilde{u}))
\frac{\v{s^\alpha(\sigma)-b^\alpha(u)}}{ l_0 + k^\alpha(u,\sigma)}
\frac{\v{s^\alpha(\sigma)-b^\alpha(\tilde{u})}}{ l_0 + k^\alpha(\tilde{u},\sigma)} \\
& \qquad \quad +B_{\tau}^\alpha(u, \tilde{u}, \sigma) \Bigm)
J_{\partial B}^\alpha(u)J_{\partial B}^\alpha(\tilde{u})J_{\partial D^\alpha}(\sigma)\,du d\tilde{u} d\sigma, 
\end{align*}
where $J_{\partial D^\alpha}(\sigma)=\sqrt{1+(g^\alpha_{\sigma_1})^2+(g^\alpha_{\sigma_2})^2}$  
and $B_{\tau}^\alpha=O(\tau^{-1})$.    
For this expression of $K_{\tau,1}^{00,\alpha}(x_0^\alpha, y_0^\alpha)$, we use Lemma \ref{a usual Laplace method}.  
Considering the definitions of $a_0$ and $b_0$ and the facts: 
\[\v{x_0^{\alpha}-y_0^{\alpha}}=l_0, \enskip  
\sqrt{\gamma_0}\partial_{\nu_x}\phi(x_0^{\alpha},y_0^{\alpha})=1,\enskip\frac{\nu_{x_0^{\alpha}}\cdot(x_0^{\alpha}-y_0^{\alpha})}{\v{x_0^{\alpha}-y_0^{\alpha}}}=-1,
\]
we find that $J_{\partial B}^\alpha(0)=J_{\partial D^\alpha}(0)=1$, $\tilde{L}^\alpha(0,0,0)=2l_0$ and 
\begin{align*}
\eta^{\alpha} (s^\alpha(0), b^\alpha(0))
\tilde{\eta}^{\alpha}(s^\alpha(0), b^\alpha(0))
&= \eta^{\alpha} (x_0^\alpha, y_0^\alpha)
\tilde{\eta}^{\alpha} (x_0^\alpha, y_0^\alpha) 
\\&
= \frac{1}{8\pi^2\gamma_0^{3/2}(l_0)^2}
\frac{\sqrt{\gamma_0} - \lambda_1(x^\alpha_0)}
{\sqrt{\gamma_0} + \lambda_1(x^\alpha_0)} \enskip\text{ for $\alpha\in\{n_+, n_-\}$,} \\
\eta^{d} (s^{d}(0), b^{d}(0))\tilde{\eta}^{d}(s^{d}(0), b^{d}(0))
&=\eta^{d} (x_0^{d}, y_0^{d})\tilde{\eta}^{d} (x_0^{d}, y_0^{d}) 
=\frac{1}{8\pi^2\gamma_0^{3/2}(l_0)^2}. 
\end{align*}
Then, we have 
\begin{align*}
K_{\tau,1}^{00,\alpha}(x_0^\alpha, y_0^\alpha)=&
\frac{e^{-\frac{\tau}{\sqrt{\gamma_0}}\tilde{L}^\alpha(0,0,0)}}
{\sqrt{\text{det}(\text{Hess}(\tilde{L}^\alpha)(0,0,0))}}
\left(\frac{2\pi\sqrt{\gamma_0}}{\tau} \right)^3\frac{\gamma_0}{\tau^2} \\
&\qquad \quad \Bigm\{ 
f(y^\alpha_0)^2\eta^{\alpha} (x_0^\alpha, y_0^\alpha)\tilde{\eta}^{\alpha} (x_0^\alpha, y_0^\alpha)
J_{\partial B}^\alpha(0)^2J_{\partial D^\alpha}(0) +O(\tau^{-1/2}) \Bigm\} \\
=&\frac{e^{-\frac{\tau}{\sqrt{\gamma_0}}2l_0}}
{2\sqrt{{\mathcal A}_\alpha({x_j^{\alpha}}; \partial{D}, \partial{B})}}
\frac{\pi\gamma_0}{\tau^5}\Big((f(y_0^\alpha))^2\frac{\sqrt{\gamma_0} - \lambda_1(x^\alpha_0)}{\sqrt{\gamma_0} + \lambda_1({x^\alpha_0})}+O(\tau^{-1/2})\Big) 
\end{align*}
for $\alpha \in \{n_+, n_-\}$, and 
\begin{align*}
K_{\tau,1}^{00,d}(x_0^\alpha, y_0^\alpha)=&
\frac{e^{-\frac{\tau}{\sqrt{\gamma_0}}\tilde{L}^{d}(0,0,0)}}
{\sqrt{\text{det}(\text{Hess}(\tilde{L}^{d})(0,0,0))}}
\left(\frac{2\pi\sqrt{\gamma_0}}{\tau} \right)^3\frac{\gamma_0}{\tau^2} \\
&\qquad \quad \Bigm\{ 
f(y_0^d)^2\eta^{d} (x_0^{d}, y_0^{d})\tilde{\eta}^{d} (x_0^{d}, y_0^{d})
J_{\partial B}^{d}(0)^2J_{\partial D^{d}}(0) +O(\tau^{-1/2}) \Bigm\} \\
=&\frac{e^{-\frac{\tau}{\sqrt{\gamma_0}}2 l_0}}
{2\sqrt{{\mathcal A}_\alpha({x_j^{\alpha}}; \partial{D}, \partial{B})}}
\frac{\pi\gamma_0}{\tau^5}((f(y_0^{d}))^2+O(\tau^{-1/2})).
\end{align*}

\par
Because of (\ref{difference}), there exists 
$A_{\tau}^{3,\alpha}(\sigma, u,\tilde{u})$ satisfying 
\begin{align*}
\tilde{K}_{\tau,1}^{0,\alpha}(x_0^\alpha, y_0^\alpha)&=\int_{\partial D{\cap}B_{r_0}(x_0^\alpha)} \int_{\tilde{U}_{y_0^\alpha}\times \tilde{U}_{y_0^\alpha}}
e^{-\frac{\tau}{\sqrt{\gamma_0}}(\v{s^\alpha(\sigma)-b^\alpha(u)}+\v{s^\alpha(\sigma)-b^\alpha(\tilde{u})})}
A_{\tau}^{3,\alpha}(\sigma, u,\tilde{u})\, du d\tilde{u} dS_x
\end{align*}
and we have $A_{\tau}^{3,\alpha}
=O(\tau^{-3})$ in the same way as 
(\ref{I_2}) and (\ref{I_3}).  
Since we have 
\begin{align*}
\v{\tilde{K}_{\tau,1}^{0,\alpha}(x_0^\alpha, y_0^\alpha)}
&\le \int_{\tilde{U}_{x_0^\alpha}}\int_{\tilde{U}_{y_0^\alpha}\times \tilde{U}_{y_0^\alpha}}
e^{-\frac{\tau}{\sqrt{\gamma_0}}\tilde{L}^\alpha(\sigma, u, \tilde{u})}\text{sup}_{(\sigma, u, \tilde{u})\in 
\tilde{U}_{x_0^\alpha}\times \tilde{U}_{y_0^\alpha}\times \tilde{U}_{y_0^\alpha}} \v{A_{\tau}^{3,\alpha}(\sigma, u,\tilde{u})}d\sigma dud\tilde{u}\\
&\le C \tau^{-3} \int_{\tilde{U}_{x_0^\alpha}\times \tilde{U}_{y_0^\alpha}\times \tilde{U}_{y_0^\alpha}}
e^{-\frac{\tau}{\sqrt{\gamma_0}}\tilde{L}^\alpha(\sigma, u, \tilde{u})}d\sigma dud\tilde{u}, 
\end{align*}
from Lemma \ref{a usual Laplace method} it follows that 
$$\v{\tilde{K}_{\tau,1}^{0,\alpha}(x_0^\alpha, y_0^\alpha)} 
\le C\tau^{-6} e^{-\frac{\tau}{\sqrt{\gamma_0}}2l_0}.$$  
From the above, we can see that the order of the top terms of $K_{\tau, 1}^\alpha$ is $e^{-\frac{\tau}{\sqrt{\gamma_0}}2l_0}\tau^{-5}$ and the sum for $\alpha$ of the top terms  
of $K_{\tau, 1}^\alpha$ is the same as the top term of $I_{\tau}$ in Theorem \ref{a result for the case l_0^+ = l_0^-_general} for $T >2l_0/\sqrt{\gamma_0}$.  
\par
Next, we estimate $K_{\tau, -k}^{\alpha} $ ($k \ge 0$).  
Note that there exists a constant $C_k> 0$ such that 
\begin{align*}
\v{K_{\tau, -k}^{\alpha}} \leq C_k \int_{\partial{D}^\alpha{\times}B{\times}B}e^{-\tau L(x,y,\tilde{y})/\sqrt{\gamma_0}}dxdyd\tilde{y}.  
\end{align*}
In the above integral, by the same argument as deriving the asymptotic form $K_{\tau, 1}^\alpha$, we obtain the following:
\begin{align*}
\v{K_{\tau, -k}^{\alpha}} \leq C_k\tau^{-5}e^{-\frac{\tau}{\sqrt{\gamma_0}}2 l_0}.
\end{align*}
\par
The lower order term than $e^{-\frac{\tau}{\sqrt{\gamma_0}}2l_0}\tau^{-4}$ can be regarded as the remainder terms. Therefore, if we admit Lemma \ref{estimates of the remainder terms for the indicator function}, it is sufficient to construct the asymptotic solution as $N=0$, and it is enough that $(\partial D, \lambda)$ is in $0$-class, which yields Theorem \ref{a result for the case l_0^+ = l_0^-_general}.  
\hfill$\square$
\par
\begin{Remark}\label{Estimates by using Laplace method}
By the same argument as above, 
for any continuous function $A(x, y, \tilde{y}; \tau)$ in $\partial{D}\times\overline{B}\times\overline{B}$ with $\tau \geq 1$, we have 
\begin{align*}
\Big\vert \int_{B\times{B}}dyd\tilde{y}\int_{\partial{D}}
&e^{-\tau\frac{L(x, y, \tilde{y})}{\sqrt{\gamma_0}}}A(x, y, \tilde{y}; \tau)dS_x \Big\vert 
\\&\leq C\sup_{(x, y, \tilde{y}) \in \partial{D}\times\overline{B}\times\overline{B}}\v{A(x, y, \tilde{y}; \tau)}\tau^{-5}e^{-\frac{2\tau}{\sqrt{\gamma_0}}l_0} \qquad(\tau \geq 1).
\end{align*}
This estimate will be used to show the remainder estimates in the next section.
\end{Remark}
\setcounter{equation}{0}
\section{Estimates of the remainder terms} 
\label{Estimates of the remainder terms}

In this section, we give a proof of Lemma \ref{estimates of the remainder terms for the indicator function} under the assumption that $(\partial D, \lambda)$ is in $N$-class for the integer $N \ge0$.  \par
We put $\tilde\Omega = \Omega$ or $\tilde\Omega = \R^3\setminus\overline{\Omega}$, and 
$$
\V{\varphi}_{H^1_\tau(\tilde\Omega)} = \big\{\V{\nabla_x\varphi}_{L^2(\tilde\Omega)}^2 
+ \tau^2\V{\varphi}_{L^2(\tilde\Omega)}^2\big\}^{1/2}. 
$$
The usual trace estimates give
\begin{align*}
\V{\varphi}_{H^{1/2}(\partial\Omega)}
\leq C\V{\varphi}_{H^1(\tilde\Omega)}
\leq C\V{\varphi}_{H^1_\tau(\tilde\Omega)}
\quad(\varphi \in H^1_\tau(\tilde\Omega), \tau \geq 1).
\end{align*}
For $L^2$-norm on $\partial\Omega$, 
we need more accurate inequalities.
\begin{Lemma}\label{trace estimates with large parameter}
Assume that $\partial\Omega$ is $C^2$. Then, there exists a constants $C > 0$ such that 
\begin{align*}
&\V{\varphi}_{L^2(\partial{\tilde\Omega})} 
\leq C\tau^{-1/2}\V{\varphi}_{H^1_\tau(\tilde\Omega)}, 
\\
&\V{\varphi}_{H^{1/2}(\partial{\tilde\Omega})} \le 
\V{\varphi}_{H^{1/2}_\tau(\partial{\tilde\Omega})} \leq 
C\V{\varphi}_{H^1_\tau(\tilde\Omega)}
\quad(\varphi \in H^1_\tau(\tilde\Omega), \tau \geq 1),  
\end{align*}
where 
$$
\V{g}_{H^{1/2}_{\tau}(\partial\tilde\Omega)} = \big\{\V{g}_{H^{1/2}(\partial\tilde\Omega)}^2
+\tau\V{g}_{L^2(\partial\tilde\Omega)}^2\big\}^{1/2}. 
$$
\end{Lemma}
A proof of Lemma \ref{trace estimates with large parameter} is given in Prop.2.1 and Appendix in \cite{M. and W. Kawashita separated}. 
\par
For $u, \varphi \in H^1(\Omega)$, we set 
\begin{align*}
{B_{\tau}}[u, \varphi] = \int_{\Omega}\big\{\gamma_0\nabla_xu(x)\cdot\nabla_x\varphi(x)
+\tau^2u(x)\varphi(x)\big\}dx 
\end{align*}
and
$$
\tilde{B}_\tau[u, \varphi] = B_{\tau}[u, \varphi] + \int_{\partial{D}^{n}}\lambda(x;\tau)u(x)\varphi(x)dS_x.
$$
Recall that for $F_{\tau, N}$, $G^n_{\tau, N}$ and $H^d_\tau$ given by (\ref{definition of F_N}), 
(\ref{definition of G_N}) and (\ref{definition of H_tau}), $\Psi_{\tau, N}^r(x, y)$ is the weak solution  of 
(\ref{equations of the remainder term (for cavities only)}) if and only if 
$\Psi_{\tau, N}^r(\cdot, y) 
\in H^1(\Omega)$ 
and it satisfies
\begin{align*}
\left\{\hskip-3pt
\begin{array}{ll}
\tilde{B}_{\tau}[\Psi_{\tau, N}^r(\cdot, y), \varphi] 
= -\langle F_{\tau, N}(\cdot, y), \varphi\rangle_{\Omega} 
- \langle G^n_{\tau, N}(\cdot, y), \varphi\rangle_{\partial{D^n}}
&(\varphi \in H^1_{0,\partial D^d}(\Omega)),
\\
\Psi_{\tau, N}^r(x, y) = H^d_\tau(x, y) &\text{on } \partial{D}^d,
\end{array}
\right.
\end{align*}
where 
\begin{align*}
\langle f, \varphi\rangle_{\Omega} = \int_{\Omega}f(x)\varphi(x)dx, \qquad
\langle g, \varphi\rangle_{\partial{D^\alpha}} = \int_{\partial{D^\alpha}}g(x)\varphi(x)dS_x \quad(\alpha = n, d).
\end{align*}
Thus, $w_{N}^r(x; \tau)$ is the weak solution satisfying 
\begin{align*}
\left\{\hskip-3pt
\begin{array}{ll}
\tilde{B}_{\tau}[w_{N}^r(\cdot; \tau), \varphi] 
= -\langle f_{\tau, N}, \varphi\rangle_{\Omega}
- \langle g^n_{\tau, N}, \varphi\rangle_{\partial{D}^n}
&\quad(\varphi \in H^1_{0,\partial D^d}(\Omega)), 
\\
w_{N}(x; \tau) = h^d_{\tau}(x) & \quad\text{on } \partial{D}^d,
\end{array}
\right.
\end{align*}
where
\begin{align*}
f_{\tau, N}(x) &= \int_{B}F_{\tau, N}(x, y)f(y)dy, 
\qquad
g^n_{\tau, N}(x) = \int_{B}G^n_{\tau, N}(x, y)f(y)dy, 
\\
h^d_{\tau}(x) &= \int_{B}H^d_{\tau}(x, y)f(y)dy.
\end{align*}
From Lemmas \ref{trace estimates with large parameter} and \ref{basic estimate for weak solutions 2}, it also follows that 
\begin{align}
\V{w_{N}^r}_{H^{1}_\tau(\Omega)} 
&\leq C\tau^{-1}
\{\V{f_{\tau, N}}_{L^2(\Omega)}
+\sqrt{\tau}\V{g^n_{\tau, N}}_{L^2(\partial D^n)}
+\tau\V{ h^d_{\tau} }_{H_{\tau}^{1/2}(\partial D^d)}\}, 
\label{estimate of V{w_{0, N}}_{H^{1}_tau(Omega)}}
\\
\V{w_{N}^r}_{L^2(\partial{D}^n)} 
&\leq C\tau^{-3/2}
\{\V{f_{\tau, N}}_{L^2(\Omega)}
+\sqrt{\tau}\V{g^n_{\tau, N}}_{L^2(\partial D^n)} 
+\tau\V{ h^d_{\tau} }_{H_{\tau}^{1/2}(\partial D^d)}\},
\label{estimate of V{w_{0, N}}_{L^2(partial D^n)}}
\end{align}
as $\tau \geq \tau_0$ for some fixed constant $\tau_0 \geq 1$ and $C > 0$.

\subsection{Estimate of $g_{\tau, N}$}

From (\ref{definition of G_N}), $G_{\tau, N}^n(x, y)$ is decomposed as
$G^n_{\tau, N}(x, y) = G_{\tau, N, 0}^{n}(x, y) + G_{\tau, N, -\infty}^{n}(x, y)$, where  
\begin{align*}
G_{\tau, N, 0}^n(x, y) &= -(-\tau)^{-N}e^{-{\tau}\phi(x, y)}\chi(x)
(\delta_{N,0}a_1(x,y)+{\mathcal B_{0}^n}b_{N}(x, y)), 
\\
G_{\tau, N, -\infty}^{n}(x, y) & = (\chi(x)-1){\mathcal B_{\tau}^n}\Phi_{\tau}(x, y) -  \gamma_0\partial_{\nu_x}\chi(x)\Psi_{\tau, N}(x, y).
\end{align*}
Since $\phi(x, y) > l_0/\sqrt{\gamma_0}$ for $(x, y) \in (({\mathcal U}_3\setminus{\mathcal U}_1)\cap\Omega)\times\overline{B}$, 
there exists a constant $c_1 >0$ such that $c_1 \le c_0$ and 
\begin{equation}
\sqrt{\gamma_0}\phi(x, y) \geq l_0 + c_1 \quad((x, y) \in (({\mathcal U}_3\setminus{\mathcal U}_1)\cap\Omega)\times\overline{B}). 
\label{estimate of phi^pm in longer regions}
\end{equation}
From (\ref{estimate of phi^pm in longer regions}) and (\ref{estimates on the outside of E^alpha}) it follows that 
$$
\v{G^n_{\tau, N, -\infty}(x, y)} \leq C_Ne^{-\frac{\tau}{\sqrt{\gamma_0}}(l_0+c_1)}
\quad((x, y) \in \partial{D}\times\overline{B})
$$
for some constant $C_N > 0$. 
Since we have
$$
\v{G^n_{\tau, N, 0}(x, y)} \leq C\tau^{-N}\chi(x)e^{-\frac{\tau}{\sqrt{\gamma_0}}\v{x - y}}\quad((x, y) \in \partial{D}\times\overline{B}),
$$
it follows that
\begin{align*}
\V{g^n_{\tau, N}}_{L^2(\partial{D})}^2 &= \int_{B{\times}B}dyd\tilde{y}f(y)f(\tilde{y})\int_{\partial{D}}G^n_{\tau, N}(x, y)G^n_{\tau, N}(x, \tilde{y})dS_x
\\&
\leq C\int_{B{\times}B}dyd\tilde{y}\int_{\partial{D}}\Big\{(\chi(x))^2e^{-\frac{\tau}{\sqrt{\gamma_0}}L(x, y, \tilde{y})}\tau^{-2N} + 3e^{-\frac{\tau}{\sqrt{\gamma_0}}(2l_0+c_1)}\Big\}dS_x,
\end{align*}
which yields 
$$
\V{g^n_{\tau, N}}_{L^2(\partial{D})}^2 \leq C\tau^{-2N - 5}e^{-\frac{2\tau}{\sqrt{\gamma_0}}l_0}
$$
by Remark \ref{Estimates by using Laplace method}.
Thus, we obtain
\begin{align}
\sqrt{\tau}\V{g^n_{\tau, N}}_{L^2(\partial{D})} \leq C\tau^{-N - 2}e^{-\frac{\tau}{\sqrt{\gamma_0}}l_0}
\quad(\tau \geq 1).
\label{estimate of sqrt{tau}V{g_{tau, N}}_{L^2(partial{D})}}
\end{align}

\subsection{Estimate of $f_{\tau, N}$}

From (\ref{definition of F_N}), $F_{\tau, N}(x, y)$ is decomposed as
${F}_{\tau, N}(x, y) = F_{\tau, N, 0}(x, y) + F_{\tau, N, -\infty}(x, y)$, where  
\begin{align*}
F_{\tau, N, 0}(x, y) &= -(-\tau)^{-N}e^{-{\tau}\phi(x, y)}\chi(x)
\gamma_0{\triangle}b_{N}(x, y), 
\\
F_{\tau, N, -\infty}(x, y) & 
= -\gamma_0[\Delta, \chi]\Psi_{\tau, N}(x, y), 
\end{align*}
where $\chi(x)\Delta b_N(x, y) \in L^\infty(\Omega\times{B})$, since $\Delta b_N(\cdot, y) \in \tilde{C}^{0}((\overline{\Omega}\cap{\mathcal U}_4)\times\overline{B})$ if $(\partial D, \lambda)$ is in $N$-class.  Choose $\tilde\chi \in C^\infty_0({\mathcal U}_3)$ such that $0 \leq \tilde\chi \leq 1$ and $\tilde\chi = 1$ in ${\mathcal U}_2$, then $\tilde\chi = 1$ near ${\rm supp}\,\chi \subset {\mathcal U}_2$.  
%
For $F_{\tau, N, -\infty}$, from (\ref{estimate of phi^pm in longer regions}) it follows that 
\begin{equation}
\v{F_{\tau, N, \varepsilon, -\infty}(x, y)} \leq C_N\tilde\chi(x)e^{-\frac{\tau}{\sqrt{\gamma_0}}(l_0+c_1)} 
\quad((x, y) \in \Omega\times\overline{B}) 
\label{part for longest region}
\end{equation}
and $F_{\tau, N, 0}$ satisfies following:
\begin{equation}
\v{F_{\tau, N, 0}(x, y)} \leq C\tau^{-N}\chi(x)e^{-{\tau}\phi(x, y)}\quad((x, y) \in \Omega\times\overline{B}).  
\label{part for near region}
\end{equation}

To estimate $f_{\tau, N}$, we begin with the estimate of $\phi$. 
From eikonal equation (\ref{eikonal eq}), we deduce 
$$
\partial_{\nu_x}\phi(x, y) = -\frac{1}{\sqrt{\gamma_0}}\frac{x - y}{\v{x - y}}\cdot\nu_x\quad((x,y) \in (\partial{D}{\cap}{\mathcal U}_3)\times \overline{B}).  
$$
Take any $(x_0^{\alpha}, y_0^{\alpha}) \in E_0^{\alpha}$.  
Since $\overline{B} \subset \{y \in \R^3 \mid (y - x_0^\alpha)\cdot\nu_{x_0^\alpha} \geq l_0 \}$, we can choose $c_2 > 0$ satisfying  
$$
\nu_{x_0^\alpha}\cdot\nabla_x\phi(x_0^\alpha, y) = \frac{1}{\sqrt{\gamma_0}}\frac{y - x_0^\alpha}{\v{x_0^\alpha - y}}\cdot\nu_{x_0^\alpha} \geq 3c_2
\quad(y \in \overline{B}).   
$$
Since $\phi \in \tilde{C}^{2N+4}((\overline{\Omega}\cap{\mathcal U}_4)\times \tilde{B})$, we have $\partial_{\nu_x}\phi(x, y) \geq 2c_2$ for $(x,y) \in (\overline{\Omega}\cap{\mathcal U}_4)\times \overline{B}$ if we take $r_0>0$ small enough.  
Also, since $\v{\nu_{x_0} - \nu_x}^2=2(1-\nu_{x_0^\alpha}\cdot\nu_{x}) \to 0$ ($x \to x_0^\alpha$), if we choose $r_0 > 0$ small again if necessary, it is valid 
\begin{equation*}
\nu_{x_0^\alpha}\cdot\nabla_x\phi(x, y) = \partial_{\nu_x}\phi(x, y) + (\nu_{x_0} - \nu_x)\cdot\nabla_x\phi(x, y) \geq c_2
\enskip((x,y) \in (\partial{D}\cap B_{3r_0}(x_0^\alpha)) \times \overline{B}).  
\end{equation*}
From properties of the characteristic equation of $\phi$, it follows 
$$
\nabla_x\phi(x + t\nabla_x\phi(x, y), y) = \nabla_x\phi(x, y) \quad(t \geq 0, (x,y) \in (\partial{D}\cap B_{3r_0}(x_0^\alpha)) \times \overline{B})$$
(cf. Chap. 3 of Ikawa \cite{{Ikawa book(asymptotic solutions) }}).
Thus we have
\begin{align*}
\nu_{x_0^\alpha}\cdot\nabla_x\phi(x, y) \geq c_2
\quad((x, y) \in (B_{3r_0}(x_0^\alpha)\cap\Omega)\times\overline{B}).  
\end{align*}
Let us denote $x \in \Omega{\cap}B_{4r_0}(x_0^\alpha)$ by the local coordinate 
$$
x = x^\alpha(\sigma, \sigma_3) = s^\alpha(\sigma) + \sigma_3\nu_{x_0^\alpha}, \quad(\sigma \in U_{x_0^\alpha}, \sigma_3 \geq 0) 
$$
and set $t_0^{\alpha}(\sigma)=\sup\{ \sigma_3 \mid x^\alpha(\sigma,\sigma_3) \in \Omega{\cap}B_{3r_0}(x_0^\alpha)\}$ $(\sigma \in \tilde{U}_{x_0^\alpha})$, then 
\begin{align}
\phi(x^\alpha(\sigma, \sigma_3), y) &= \phi(s^\alpha(\sigma), y) + \int_0^{\sigma_3}\nu_{x_0^\alpha}\cdot\nabla_x\phi(x^\alpha(\sigma, t), y)dt
\nonumber
\\&
\geq \frac{\v{s^\alpha(\sigma) - y}}{\sqrt{\gamma_0}} + c_2\sigma_3
\quad(\sigma \in \tilde{U}_{x_0^\alpha}, 0 \leq \sigma_3 \leq t^\alpha_0(\sigma)).
\label{est of phi nearest part}
\end{align}
Moreover, there exists a constant $c_1' > 0$ such that
\begin{equation}
\sqrt{\gamma_0}\phi(x, y) \geq l_0 + c_1' 
\quad((x, y) \in (\Omega \cap B_{3r_0}(x_0^\alpha))\times(\overline{B}\setminus{B_{r_0}(y_0^\alpha)})).
\label{estimate of phi^pm in longer regions 2}
\end{equation}
We show the following estimate:
\begin{Lemma}\label{estimate for f tauN}
For the solution $\phi$ of (\ref{eikonal eq}), there exists a constant $C>0$ such that 
\begin{equation}\int_{{B}\times{B}}dyd\tilde{y} \int_{\Omega}
(\chi(x))^2e^{-\tau(\phi(x,y)+\phi(x,\tilde{y}))}dx \le 
C\tau^{-6}e^{-\frac{2\tau}{\sqrt{\gamma_0}}l_0} \quad
(\tau \ge 1). \label{est for f tauN}
\end{equation}
\end{Lemma}
{\it Proof of Lemma \ref{estimate for f tauN}.}
From (\ref{estimate of phi^pm in longer regions 2}),
it follows that the left-hand side of (\ref{est for f tauN}) can be evaluated from above by 
\begin{align}
 \sum_{\alpha \in \{n^+, n^-,d\}} \int_{(B{\cap}B_{r_0}(y_0^\alpha))^2}dyd\tilde{y}
 \int_{\Omega\cap B_{r_0}(x_0^{\alpha})}
 (\chi(x))^2e^{-\tau(\phi(x,y)+\phi(x,\tilde{y}))}dx
+\tilde{C} e^{-\frac{\tau}{\sqrt{\gamma_0}}(2l_0+c_1')}. 
\label{1st est for f tauN}
\end{align}
Here, from (\ref{est of phi nearest part}) for $y, \tilde{y} \in B$ we have 
\begin{align*}
\int_{\Omega\cap B_{r_0}(x_0^{\alpha})}(\chi(x))^2e^{-{\tau}(\phi(x, y) + \phi(x, \tilde{y}))}dx \leq \int_{\tilde{U}_{x_0}}\int_{0}^{t^\alpha_0(\sigma)}e^{-2c_2{\tau}\sigma_3}e^{-\frac{\tau}{\sqrt{\gamma_0}}(\v{s^\alpha(\sigma) - y}+\v{s^\alpha(\sigma) - \tilde{y}})}d\sigma_3d\sigma.  
\end{align*}
Thus, the integral in (\ref{1st est for f tauN}) is estimated by Remark \ref{Estimates by using Laplace method} as 
\begin{align*}
\frac{1}{2c_2\tau}\int_{(B{\cap}B_{r_0}(y_0^\alpha))^2}\int_{\tilde{U}_{x_0}}e^{-\frac{\tau}{\sqrt{\gamma_0}}(\v{s^\alpha(\sigma) - y}+\v{s^\alpha(\sigma) - \tilde{y}})}d{\sigma}dyd\tilde{y}
\leq C\tau^{-6}e^{-\frac{2\tau}{\sqrt{\gamma_0}}l_0}.  
\end{align*}
Thus, we obtain Lemma \ref{estimate for f tauN}.  
\hfill$\square$
\par
As well as $g_{\tau, N}$, it follows from (\ref{part for longest region}), (\ref{part for near region}) and (\ref{estimate of phi^pm in longer regions 2}) that
\begin{align*}
\V{&f_{\tau, N}}_{L^2(\Omega)}^2 = \int_{B{\times}B}dyd\tilde{y}f(y)f(\tilde{y})\int_{\Omega}F_{\tau, N}(x, y)F_{\tau, N}(x, \tilde{y})dx
\\&
\leq C\int_{B{\times}B}dyd\tilde{y}\int_{\Omega}\Big\{(\chi(x))^2e^{-{\tau}(\phi(x, y) + \phi(x, \tilde{y}))}\tau^{-2N} 
+ 3e^{-\frac{\tau}{\sqrt{\gamma_0}}(2l_0+c_1)}\Big\}dx.
\end{align*}
From this estimate and Lemma \ref{estimate for f tauN},  
we obtain 
\begin{align}
\V{f_{\tau, N}}_{L^2(\Omega)} \leq C\tau^{-N - 3}e^{-\frac{\tau}{\sqrt{\gamma_0}}l_0}
\quad(\tau \geq 1).
\label{estimate of V{f_{tau, N}}_{L^2(Omega)} for C^{N+4}}
\end{align}

\subsection{Estimates of $h_{\tau}$}

From (\ref{definition of H_tau}), considering $\chi(x)-1=0$ in $\mathcal{U}_1$ and (\ref{estimates on the outside of E^alpha}), we have 
\begin{align*}
\V{h^d_{\tau}}_{L^2({D^d})}^2 
&= \int_{B{\times}B}dyd\tilde{y}f(y)f(\tilde{y})\int_{D^d}(\chi(x)-1)^2\Phi_\tau(x, y)\Phi_\tau(x, \tilde{y})dx 
\le C e^{-\frac{2\tau}{\sqrt{\gamma_0}}(l_0+c_0)}. 
\end{align*}
Since we have
$\partial_{x_j} \Phi_{\tau} =
\frac{-\tau e^{-\tau\v{x-y}/\sqrt{\gamma_0}}}{4\pi\gamma_0\v{x-y}^2}
\Big(\frac{1}{\sqrt{\gamma_0}}+\frac{1}{\tau \v{x-y}}
\Big)(x_j-y_j)$, 
in the same way as above, it follows 
$$\V{\nabla h^d_{\tau}}_{L^2({D^d})}^2  \le \tilde{C}
{\tau}^2 e^{-\frac{2\tau}{\sqrt{\gamma_0}}(l_0+c_0)}.  $$
Thus, Lemma \ref{trace estimates with large parameter} 
implies that 
\begin{equation}
\V{h^d_{\tau}}_{H_{\tau}^{1/2}(\partial D^d)} 
\le C\V{h_{\tau}}_{H_{\tau}^{1}(D^d)}  
\le C \tau e^{-\frac{\tau}{\sqrt{\gamma_0}}(l_0+c_0)}.
\label{estimate of h_tau}
\end{equation}

\subsection{Estimates of $J_{\tau, N}^{\alpha,r}$}

We begin with the estimate of $v$.  
\begin{Lemma}\label{estimates of v on partial{D}}
If $B$ is a convex set with $C^2$ boundary and $D$ and $B$ satisfy the non-degenerate condition, then there are the following estimates:   
\begin{align}
\V{v(\cdot; \tau)}_{L^{2}(\partial{D})} &\leq C\tau^{-5/2}e^{-\frac{\tau}{\sqrt{\gamma_0}}l_0} \quad(\tau \geq 1), \label{est of v}
\\
\V{\partial_{\nu_x}v(\cdot; \tau)}_{L^2(\partial{D})} 
&\leq C\tau^{-3/2}e^{-\frac{\tau}{\sqrt{\gamma_0}}l_0} \quad(\tau \geq 1), 
\label{est of partial nu v}
\\
\V{v(\cdot; \tau)}_{H^{1/2}_\tau(\partial{D})} &\leq C\tau^{-2}e^{-\frac{\tau}{\sqrt{\gamma_0}}l_0} \quad(\tau \geq 1). 
\label{est of partial v on boundary}
\end{align}
\end{Lemma}

{\it Proof.} From (\ref{free kernel_rep}) it follows that
\begin{align*}
\V{v(\cdot; \tau)}_{L^2(\partial{D})}^2 &
= \int_{B{\times}B}dyd\tilde{y}f(y)f(\tilde{y})\int_{\partial{D}}e^{-\frac{\tau}{\sqrt{\gamma_0}}L(x, y, \tilde{y})}A(x, y, \tilde{y})dS_x,
\end{align*}
where
$
A(x, y, \tilde{y}) = \frac{1}{16\pi^2\gamma_0^2}\frac{1}{\v{x - y}\v{x - \tilde{y}}}  
$.  
By Remark \ref{Estimates by using Laplace method}, we obtain (\ref{est of v}).  In the same way, from 
$$
\partial_{\nu_x}v(x; \tau) = -\tau\int_{B}e^{-\frac{\tau}{\sqrt{\gamma_0}}\v{x - y}}A_1(x, y; \tau)f(y)dy,
$$
where
$
A_1(x, y; \tau) = \frac{1}{4\pi\gamma_0^{3/2}}\Big(1 + \frac{\sqrt{\gamma_0}}{\tau}\frac{1}{\v{x-y}}\Big)\frac{\nu_x\cdot(x - y)}{\v{x-y}^2}
$, 
it follows that
\begin{align*}
\V{\partial_{\nu_x}v(\cdot; \tau)}_{L^2(\partial{D})}^2 = \tau^2\int_{B{\times}B}dyd\tilde{y}f(y)f(\tilde{y})\int_{\partial{D}}e^{-\frac{\tau}{\sqrt{\gamma_0}}L(x, y, \tilde{y})}A_1(x, y; \tau)A_1(x, \tilde{y}; \tau)dS_x.  
\end{align*}
Thus, by Remark \ref{Estimates by using Laplace method} we have (\ref{est of partial nu v}).  
Since Lemma \ref{trace estimates with large parameter} implies that 
\begin{align*}
\V{v(\cdot; \tau)}_{H^{1/2}_\tau(\partial{D})}^2 
&\leq C\V{v(\cdot; \tau)}_{H^1_\tau(D)}^2 
= C\{\V{\nabla_xv(\cdot; \tau)}_{L^2(D)}^2 + \tau^2\V{v(\cdot; \tau)}_{L^2(D)}^2\}, 
\end{align*}
if we note $(\gamma_0\triangle - \tau^2)v = 0$ in $D$, we have 
\begin{align*}
\gamma_0\V{\nabla_xv(\cdot; \tau)}_{L^2(D)}^2 &+ \tau^2\V{v(\cdot; \tau)}_{L^2(D)}^2
= \int_{D}\big(\gamma_0\nabla_xv\cdot\nabla_xv + \tau^2\v{v}^2\big)dx
\\&
= \int_{\partial{D}}\gamma_0\partial_{\nu_x}v\cdot{v}dS_x + \int_{D}(-\gamma_0\triangle{v} + \tau^2v)vdx
\\&
\leq C\V{\partial_{\nu_x}v(\cdot ; \tau)}_{L^2(\partial{D})}\V{v(\cdot ;\tau)}_{L^2(\partial{D})} \leq C\tau^{-4}e^{-\frac{\tau}{\sqrt{\gamma_0}}2l_0}, 
\end{align*}
which implies that (\ref{est of partial v on boundary}).  %
\hfill$\square$
\par
\vskip1pc
Now, we are ready to prove Lemma \ref{estimates of the remainder terms for the indicator function}.  
%
Note that from (\ref{estimate of V{w_{0, N}}_{L^2(partial D^n)}}), 
(\ref{estimate of sqrt{tau}V{g_{tau, N}}_{L^2(partial{D})}}), (\ref{estimate of V{f_{tau, N}}_{L^2(Omega)} for C^{N+4}}) and (\ref{estimate of h_tau}) it follows that 
\begin{align}
\V{w_{N}^r}_{L^2(\partial{D}^n)} 
&\leq C\tau^{-N-7/2}e^{-\frac{\tau}{\sqrt{\gamma_0}}l_0}.  
\quad(\tau \geq \tau_0).  \label{est of w_N^r}
\end{align}
%
%
%
First, consider $J^{n_\pm,r}_{\tau, N}$.  Since 
$\V{\lambda(\cdot ; \tau)}_{L^\infty(\partial{D}^n)} \leq C\tau$ $(\tau \geq 1)$, by (\ref{est of w_N^r}) and Lemma \ref{estimates of v on partial{D}} we have 
\begin{align*}
\sum_{\alpha \in \{n^+, n^-\}}\v{J^{\alpha,r}_{\tau, N}} &\leq \sum_{\alpha \in \{n^+, n^-\}}\Big\vert\int_{\partial{D^\alpha}}{\mathcal B_{\tau}^n}v(x; \tau) 
w_{N}^r(x; \tau)dS_x\Big\vert
\\&
\leq C(\V{\partial_{\nu_x}v(\cdot; \tau)}_{L^2(\partial{D}^n)} + \tau\V{v(\cdot ; \tau)}_{L^2(\partial{D}^n)})\V{w_{N}^r(\cdot; \tau)}_{L^2(\partial{D}^n)}
\\&
\leq C(\tau^{-3/2} + \tau\cdot\tau^{-5/2})\tau^{-N-7/2}e^{-2\frac{\tau}{\sqrt{\gamma_0}}l_0}
\leq C\tau^{-N-5}e^{-2\frac{\tau}{\sqrt{\gamma_0}}l_0}.  
\end{align*}
For the estimate of $J^{d,r}_{\tau, N}$, we consider 
$\gamma_0\partial_{\nu_x}w_{N}^r$ on $\partial{D}^d$, which is given by 
\begin{align} 
\int_{\partial{D}^d}\gamma_0\partial_{\nu_x}w_{N}^r(x; \tau)\varphi(x)&dS_x
= -{B}_{\tau}[w_{N}^r(\cdot; \tau), \varphi] 
-\langle f_{\tau, N}, \varphi\rangle_{\Omega}
\nonumber\\
&- \langle {\lambda}(\cdot ; \tau)w_{N}^r(\cdot; \tau)+ g^n_{\tau, N}, \varphi\rangle_{\partial{D}^n}
\qquad(\varphi \in H^1(\Omega)). 
\label{definition of the conormal derivative}
\end{align}
(cf.  Lemma 4.3 of \cite{McLean}). 
From (\ref{definition of the conormal derivative}), Lemma \ref{trace estimates with large parameter} and (\ref{estimate of V{w_{0, N}}_{H^{1}_tau(Omega)}}), it follows that
\begin{align*}
\Big\vert
&\int_{\partial{D}^d}\gamma_0\partial_{\nu_x}w_{N}^r(x; \tau)\varphi(x)dS_x
\Big\vert
\leq C\V{w_{N}^r(\cdot; \tau)}_{H^1_\tau(\Omega)}\V{\varphi}_{H^1_\tau(\Omega)}
+ \V{f_{\tau, N}}_{L^2(\Omega)}\V{\varphi}_{L^2(\Omega)}
\nonumber\\
&\qquad\quad + \big\{\frac{\V{\lambda(\cdot ; \tau)}_{L^\infty(\partial{D}^n)}}{\tau}
\tau^{1/2}\V{w_{N}^r(\cdot; \tau)}_{L^2(\partial{D}^n)}
+ \tau^{-1/2}\V{g^n_{\tau, N}}_{L^2(\partial{D}^n)}\big\}
\tau^{1/2}\V{\varphi}_{L^2(\partial{D}^n)}
\\&
\leq C\tau^{-1}\{\V{f_{\tau, N}}_{L^2(\Omega)}
+ \sqrt{\tau}\V{g^n_{\tau, N}}_{L^2(\partial{D}^n)}
+\tau\V{ h^d_{\tau} }_{H_{\tau}^{1/2}(\partial D^d)}\}\V{\varphi}_{H^1_\tau(\Omega)},  
\end{align*}
which yields
\begin{align}
\Big\vert
\int_{\partial{D}^d}&\gamma_0\partial_{\nu_x}w_{N}^r(x; \tau)\varphi(x)dS_x
\Big\vert
\leq C\tau^{-N-3}e^{-{\tau}l_0/\sqrt{\gamma_0}}\V{\varphi}_{H^1_\tau(\Omega)},  
\label{estimate of the conormal derivative of w_{0, N} no1}
\end{align}
for $\varphi \in H^1(\Omega)$ and $\tau \geq \tau_0$ 
by (\ref{estimate of sqrt{tau}V{g_{tau, N}}_{L^2(partial{D})}}),  (\ref{estimate of V{f_{tau, N}}_{L^2(Omega)} for C^{N+4}}) and 
(\ref{estimate of h_tau}).
\par
\vskip1pc
Here, we recall an extension of functions on $\partial{D}^d$ to $\Omega$.
\begin{Prop}\label{extended opt with large parameter}
Assume that $\partial{D}^d$ is $C^2$. Then there exists a linear extension operator $E_\tau: H^{1/2}(\partial{D}^d) \to H^{1}(\Omega)$ satisfying
$E_{\tau}g = g$ on $\partial{D}^d$ in the trace sense and
$$
\V{E_{\tau}g}_{H^1_\tau(\Omega)} \leq C\V{g}_{H^{1/2}_{\tau}(\partial{D}^d)}
\qquad(\tau \geq 1, g \in H^{1/2}(\partial{D}^d)).
$$
\end{Prop}
(cf. (2) of Proposition 2.1 in \cite{M. and W. Kawashita separated})

\par

Take any $g \in H^{1/2}(\partial{D}^d)$ and put $\varphi = E_{\tau}g \in H^1(\Omega)$ in (\ref{estimate of the conormal derivative  of w_{0, N} no1}). Then, from Proposition \ref{extended opt with large parameter}, it follows that
\begin{align}
\Big\vert
\int_{\partial{D}^d}\gamma_0\partial_{\nu_x}w_{N}^r(x; \tau)g(x)dS_x
\Big\vert
\leq C\tau^{-N-3}&e^{-{\tau}l_0/\sqrt{\gamma_0}}\V{g}_{H^{1/2}_{\tau}(\partial{D}^d)} 
\label{dual forms estimate of the conormal derivative of w_{0, N} no2}
\end{align}
for $g \in H^{1/2}(\partial{D}^d)$. 
%
%
From (\ref{dual forms estimate of the conormal derivative of w_{0, N} no2}) and Lemma \ref{estimates of v on partial{D}}, we obtain
\begin{align*}
\v{J_{\tau, N}^{d,r}} &\leq \gamma_0\Big\vert\int_{\partial{D^{d}}}
\partial_{\nu_x}w_{N}^r(x; \tau)v(x; \tau)dS_x\Big\vert
\leq C\tau^{-N-3}e^{-{\tau}l_0/\sqrt{\gamma_0}}\V{v(\cdot; \tau)}_{H^{1/2}_\tau(\partial{D}^d)}
\\&
\leq C\tau^{-N-3}\tau^{-2}e^{-2{\tau}l_0/\sqrt{\gamma_0}}
= C\tau^{-N-5}e^{-2{\tau}l_0/\sqrt{\gamma_0}}
\quad(\tau \geq \tau_0).
\end{align*}
Thus, we obtain Lemma \ref{estimates of the remainder terms for the indicator function}.
\hfill$\square$

\setcounter{equation}{0}
\appendix
\renewcommand{\theequation}{A.\arabic{equation}}
\section{Appendix.}
\label{Appendix. }

\subsection{Estimates of weak solutions}

Here, we give some estimates for a weak solution $w \in H^1_{0, \partial{D}^d}(\Omega)$ of 
\begin{align}
\tilde{B}_\tau[w, \varphi] = - \langle f_0, \varphi\rangle_{\Omega} - \sum_{j = 1}^3\langle f_j, \partial_{x_j}\varphi\rangle_{\Omega} - \langle g, \varphi\rangle_{\partial{D}^n}
\quad (\varphi \in H^1_{0, \partial{D}^d}(\Omega)),
\label{weak formulation of elliptic problems}
\end{align}
where $f_j \in L^2(\Omega)$ $(j = 0, 1, 2, 3)$ and $g \in L^2(\partial{D}^n)$.
\begin{Lemma}\label{basic estimate for weak solutions}
There exist constants $C_0 > 0$ and $\tau_0 > 0$ such that
\begin{align*}
\V{w}_{H^1_\tau(\Omega)} &\leq C_0\tau^{-1}
\{\V{f_{0}}_{L^2(\Omega)}+ \tau\sum_{j = 1}^3\V{f_{j}}_{L^2(\Omega)}
+ \sqrt{\tau}\V{g}_{L^2(\partial{D}^n)}\}
\quad(\tau \geq \tau_0).
\end{align*}
\end{Lemma}
\par
{\it Proof.} Taking $\varphi = w \in H^1_{0,\partial D^d}(\Omega)$ in 
(\ref{weak formulation of elliptic problems}), we have
\begin{align*}
B_{\tau}[w, w] + \tau\langle \lambda_1w, w\rangle_{\partial{D}^n}
= &-\langle f_{0}, w \rangle_{\Omega} 
- \sum_{j = 1}^3\langle f_{j}, \partial_{x_j}w \rangle_{\Omega}
\\&
- \langle g, w \rangle_{\partial{D}^n} 
- \langle \lambda_0w, w \rangle_{\partial{D}^n}
\quad(\varphi \in H^1_{0, \partial{D}^d}(\Omega))
\nonumber
\end{align*}
Since $\lambda_1(x) \ge 0$, we have
\begin{align*}
\V{w}_{H^1_\tau(\Omega)}^2
\leq C_1&\big\{\V{f_0}_{L^2(\Omega)}\V{w}_{L^2(\Omega)}
+ \sum_{j = 1}^3\V{f_j}_{L^2(\Omega)}\V{\partial_{x_j}w}_{L^2(\Omega)}
\\&
+\V{\lambda_0}_{L^\infty(\partial{D}^n)}\V{w}_{L^2(\partial D^n)}^2
+\V{g}_{L^2(\partial D^n)}\V{w}_{L^2(\partial D^n)}\big\}, 
\end{align*}
where $C_1 = \min\{1/\gamma_0, 1\} > 0$.
From this estimate and Lemma \ref{trace estimates with large parameter}, it follows that there exists a constant $C > 0$ such that
\begin{align*}
\V{w}_{H^1_\tau(\Omega)}^2 &\leq C_1\V{w}_{H^1_\tau(\Omega)}
\Big\{{\tau}^{-1}\V{f_0}_{L^2(\Omega)}
+ \sum_{j = 1}^3\V{f_{j}}_{L^2(\Omega)}
+ C{\tau}^{-\frac{1}{2}}\V{g}_{L^2(\partial D^n)}\Big\}
\\&\hskip30mm
+ \frac{C_1C^2}{\tau}\V{\lambda_0}_{L^\infty(\partial D^n)}
\V{w}_{H^1_\tau(\Omega)}^2
\\&
\leq C_1^2\tau^{-2}\Big\{\V{f_0}_{L^2(\Omega)}
+ \tau\sum_{j = 1}^3\V{f_{j}}_{L^2(\Omega)}
+ C\sqrt{\tau}\V{g}_{L^2(\partial D^n)}\Big\}^2
\\&\hskip30mm
+ \Big(\frac{1}{4} + \frac{C_1C^2}{\tau}\V{\lambda_0}_{L^\infty(\partial D^n)}
\Big)\V{w}_{H^1_\tau(\Omega)}^2. 
\end{align*}
Hence, we obtain Lemma \ref{basic estimate for weak solutions} if we take $\tau \geq \max\{4C_1C^2\V{\lambda_0}_{L^\infty(\partial D^n)}, 1\}$.
\hfill
$\square$
\vskip1pc\par\noindent
If $w \in H^1(\Omega)$ satisfies (\ref{weak formulation of elliptic problems}) for $f_0$, $f_j \in L^2(\Omega)$, $g \in L^2(\partial{D}^n)$ and $w\vert_{\partial{D}^d} = h$ for $h \in H^{1/2}(\partial{D}^d)$, we call $w$ a weak solution for the equations: 
\begin{equation}
\left\{
\begin{array}{ll}
({\gamma_0}\Delta - \tau^2)w = f_0 -\sum_{j = 1}^3\partial_{x_j}f_j
&\qquad\text{in } \Omega, \\
{\mathcal B_{\tau}^n}w = g
& \qquad\text{on } \partial{D}^n,  
\\ 
w = h
& \qquad\text{on } \partial{D}^{d}.  
\end{array}
\right.
\label{BVP with inhomogeneous Dirichlet condition}
\end{equation}
The uniqueness of the solution $w$ for (\ref{BVP with inhomogeneous Dirichlet condition}) follows the uniqueness of solution $w \in H^1_{0, \partial{D}^d}(\Omega)$ for (\ref{weak formulation of elliptic problems}). 
To show the existence of the solution for (\ref{BVP with inhomogeneous Dirichlet condition}), we consider $\tilde{w} = w - \chi E_{\tau}h \in H^1(\Omega)$ by using the operator $E_{\tau}$ of Proposition \ref{extended opt with large parameter}, where $\chi$ is  introduced in (\ref{definition of the remainder term}).  Since $\tilde{w}\vert_{\partial{D}^d} = 0$, we have $\tilde{w} \in H^1_{0, \partial{D}^d}(\Omega)$, and $\tilde{w}$ satisfies 
\begin{align}
\tilde{B}_\tau[\tilde{w}, \varphi] 
=& - \langle f_0, \varphi\rangle_{\Omega} - \sum_{j = 1}^3\langle f_j, \partial_{x_j}\varphi\rangle_{\Omega} - \langle g, \varphi\rangle_{\partial{D}^n}
\nonumber
\\&- \sum_{j = 1}^3\langle\partial_{x_j}(\chi E_{\tau}h), \partial_{x_j}\varphi\rangle_{\Omega} 
- \tau^2\langle\chi E_{\tau}h, \varphi\rangle_{\Omega}, 
\label{weak form of A.4}
\end{align}
for any $\varphi \in H^1_{0, \partial{D}^d}(\Omega)$.  
The existence of the weak solution $\tilde{w} \in H^1_{0, \partial{D}^d}(\Omega)$ for (\ref{weak form of A.4}) implies the existence of the solution $w \in H^1(\Omega)$ for (\ref{BVP with inhomogeneous Dirichlet condition}). \par
Let us estimate the weak solution $w$ of (\ref{BVP with inhomogeneous Dirichlet condition}). 
From Lemma \ref{basic estimate for weak solutions}, for $\tau \geq \tau_0$ we have 
\begin{align*}
\V{w}_{H^1_\tau(\Omega)} &\leq C_0\tau^{-1}
\Big\{\V{f_{0}}_{L^2(\Omega)}+ \tau\sum_{j = 1}^3\V{f_{j}}_{L^2(\Omega)}
+ \sqrt{\tau}\V{g}_{L^2(\partial{D}^n)}
\\&\hskip20mm 
+ \tau\sum_{j = 1}^3\V{\partial_{x_j}(\chi E_{\tau}h)}_{L^2(\Omega)} 
+\tau^2\V{\chi E_{\tau}h}_{L^2(\Omega)}\Big\}+ \V{\chi E_{\tau}h}_{H^1_\tau(\Omega)}
\\&
\leq C_0\tau^{-1}
\Big\{\V{f_{0}}_{L^2(\Omega)}+ \tau\sum_{j = 1}^3\V{f_{j}}_{L^2(\Omega)}
+ \sqrt{\tau}\V{g}_{L^2(\partial{D}^n)}\Big\} + C\V{E_{\tau}h}_{H^1_{\tau}(\Omega)}. 
\end{align*}
From the above and Proposition \ref{extended opt with large parameter} it follows that 
\begin{Lemma}\label{basic estimate for weak solutions 2}
There exist constants $C_0 > 0$ and $\tau_0 > 0$ such that
\begin{align*}
\V{w}_{H^1_\tau(\Omega)} &\leq C_0\tau^{-1}
\{\V{f_{0}}_{L^2(\Omega)}+ \tau\sum_{j = 1}^3\V{f_{j}}_{L^2(\Omega)}
+ \sqrt{\tau}\V{g}_{L^2(\partial{D}^n)} + \tau\V{h}_{H^{1/2}_{\tau}(\partial{D}^d)}\}
\end{align*}
for $\tau \geq \tau_0$. 
\end{Lemma}

\bibliographystyle{plain}
\bibliography{bibfile}

\vskip1pc
\footnotesize
\par\noindent
{\sc Kawashita, Mishio}
\par\noindent
{\sc Mathematics Program, \\
Graduate School of Advanced Science and Engineering,
\\
Hiroshima University, Higashihiroshima 739-8526, Japan.}
\vskip1pc
\footnotesize
\par\noindent
{\sc Kawashita, Wakako}
\par\noindent
{\sc Electrical, Systems, and Control Engineering Program, \\
Graduate School of Advanced Science and Engineering, \\
Hiroshima University, 
Higashihiroshima 739-8527, Japan.}

\end{document}